A Fragile Points Method, with an interface debonding model, to simulate damage and fracture of U-notched structures

Kailei Wang[1], Baoying Shen[1], Mingjing Li[1*], Leiting Dong [1*], Satya N. Atluri[2]

[1] School of Aeronautic Science and Engineering, Beihang University, Beijing, CHINA

[2] Department of Mechanical Engineering, Texas Tech University, Lubbock, Texas, USA

**Abstract:**

Notched components are commonly used in engineering structures, where stress concentration may easily lead to crack initiation and development. The main goal of this work is to develop a simple numerical method to predict the strength and crack-growth-path of U-notched specimens made of brittle materials. For this purpose, the Fragile Points Method (FPM), as previously proposed by the authors, has been augmented by an interface debonding model at the interfaces of the FPM domains, to simulate crack initiation and development. The formulations of FPM are based on a discontinuous Galerkin weak form where point-based piece-wise-continuous polynomial test and trial functions are used instead of element-based basis functions. In this work, the numerical fluxes introduced across interior interfaces between subdomains are postulated as the tractions acting on the interface derived from an interface debonding model. The interface damage is triggered when the numerical flux reaches the interface strength, and the process of crack-surface separation is governed by the fracture energy. In this way, arbitrary crack initiation and propagation can be naturally simulated without the need for knowing the fracture-patch before-hand. Additionally, a small penalty parameter is sufficient to enforce the weak-form continuity condition before damage initiation, without causing problems such as artificial compliance and numerical ill-conditioning. As validations, the proposed FPM method with the interface debonding model is used to predict fracture strength and crack-growth trajectories of U-notched structures made of brittle materials, which is useful but challenging in engineering structural design practices.

---
* Email address of ccorresponding authors: limingjing@buaa.edu.cn (M. Li) and ltdong@buaa.edu.cn (L. Dong)

# 1 Introduction

Notches are widely-used in engineering structures, where stress concentration may cause crack-initiation and structural failure. Thus, during the past decades, damage tolerance assessment of notched structures has attracted many research efforts, most of which have been focused on the following two aspects. The first one is predicting the fracture strength of designed structures in order to determine their load bearing capacity. The second one is predicting the crack-growth trajectories in order to efficiently reinforce the designed structures or to guide the crack path in a desired direction.

For the assessment of fracture strength, linear elastic fracture mechanics (LEFM) has been applied for sharp notches. The stress field near the root can be characterized by the intensity of its singularity, and the failure criterion can be therefore postulated. However, when the notch root is blunt, i.e. the notch root radius is not zero, the stress singularity disappears and analyses using LEFM are overly conservative [1]. In this work, we focus on damage assessment of U-notched specimens made of brittle materials, for which the failure criterion based on stress intensity factor gives poor results [2]. Several empirical failure criteria have been proposed to predict the load-bearing capacities of U-notched structures [3, 4]. However, for more complex structures with various notch shapes, root radii and loading modes, predicting fracture strength of U-notched specimens by means of simple analytical formulae is challenging, and it is also more difficult when it comes to predicting the entire process from crack initiation to structural fracture.

During the past decades, several numerical approaches based on finite element method (FEM) have been developed to assess the fracture resistance of U-notched specimens. In general, these approaches can be classified into two categories.

In the first category, certain failure criteria, such as the mean stress criterion [5-7] and the strain energy density criterion [8-10], have been applied to determine failure in combination with FEM. In these approaches, failure was triggered by certain variables, e.g. stress, stress intensity factor and strain energy density, evaluated within elements near the notch root, and the maximum load was estimated when the variables reached the prescribed critical value.

The second category of approaches assume that cracks are triggered by the accumulation of microscopic damage. The most representative and widely-used approach is based on cohesive zone model (CZM) [11-13], in which cracks are modeled explicitly by interface cohesive elements and a traction-separation law is defined to characterize damage initiation and development between neighboring solid elements. CZM has been successfully applied to predict fracture strength of U-notched specimens made of PMMA at low temperature [14] and graphite [15].

Though numerical approaches reviewed above can provide good estimations on fracture strength, the prediction on fracture trajectory has been barely considered. The first category of numerical approaches reviewed above is not able to model crack propagation without additional mesh regeneration technique. CZM can model both crack initiation and propagation. However, when the crack path is unknown in advance, a large number of interface cohesive elements has to be inserted into the finite element (FE) models leading to a dramatic increase of degrees of freedom (DOFs) and leading to the numerical problem known as artificial compliance [16]. To avoid artificial compliance, a very large initial interfacial stiffness or penalty function has to be used, resulting in ill-conditioning in implicit analysis and the decrease in stable time increment in explicit analysis [17, 18].

For predicting the fracture trajectory, most existing FE-based numerical approaches are based on the incremental method and remeshing techniques [19-21].

In these approaches, the location and orientation of crack growth are determined based on certain criteria such as the maximum principal tensile stress criterion [22-24], and another criterion, normally the maximum tangential stress (MTS) criterion [25], are usually used to determine the fracture angle during crack propagation. During the simulation, remeshings are used repeatedly to accommodate inserted cracks, and singular elements are often used in order to improve the computed stress field near the crack-tip[21]. Extended finite element method (XFEM) is also proposed and used, with crack-surface and crack-tip enrichments, to alleviate the need for remeshing[26-29]. XFEM has also been combined with CZMs to predict fracture trajectory of various blunt-notched specimens [30].

Recently, a simple meshless method, named the Fragile Points Method (FPM), has been proposed by the group of Dong and Atluri [31, 32]. In this method, the problem domain is discretized by points and the corresponding formulations are based on the Discontinuous Galerkin weak form, different from other weak form meshless methods, e.g. meshless local Petrov-Galerkin (MLPG) method [33] and element-free Galerkin (EFG) method [34]. Instead of using element-based basis functions, the trial and test functions used in FPM are point-based piece-wise-continuous polynomials, for which the integral of the Galerkin weak form can be evaluated efficiently by using a very low-order Gauss quadrature. Additionally, as the FPM trial functions are discontinuous across the interfaces between neighboring subdomains, it is very convenient to introduce cracks simply by manipulating the numerical fluxes and support domains without changing the number of DOFs. In previous studies, FPM has been applied to solve Poisson's equation [31] and linear elastic problems [32] to verify the its stability, robustness, accuracy and efficiency. This method has also been used to simulate transient heat conduction in anisotropic nonhomogeneous media [35, 36] and flexoelectric effects in dielectric solids [37].

To predict fracture strength and trajectory of fracture propagation of U-notched structures made of brittle materials, FPM formulations have been enriched in this work by an interface debonding model. Note that the symmetric interior penalty numerical flux was employed in the FPM formulation by the previous works of the authors, see [31]. Instead, we use the incomplete interior penalty Galerkin (IIPG) [38] numerical flux in this work and point out the IIPG numerical flux has an explicit physical meaning as the traction acting on interior interfaces. Thus, the process of interface damage and crack initiation can easily be modeled by defining a damaged interface numerical flux postulated as the traction derived from an interface debonding model.

In this work, such an interface debonding model is postulated by following the available previous works of thermodynamically consistent Continuum Damage Models (CDMs). According to the proposed interface debonding model, damage is initiated across an interior interface when the traction reaches the material strength, and the separation process is governed by the material fracture energy. Compared to the intrinsic CZM, interior interfaces in FPM models remain closed (in a weak sense) before reaching the material strength, even though a small penalty factor is used. Thus, this method does not suffer from artificial compliance, and avoids the need for using an extremely large value of interfacial stiffness.

Note that the proposed FPM with interface debonding model shares some similarities with discontinuous Galerkin (DG) FEM in combine with CZM [39-41]. However, in hybrid DG-FEM/CZM approaches, cohesive interface elements are inserted into damaged interior interfaces, which means that, for problems with dense crack initiations, many additional nodes are inserted into the model and the number of DOFs is increased continuously. In contrast, the proposed FPM does not require inserting extra nodes, and can be easily used to simulate the entire process from small crack initiation from notch root to the structural failure caused by the propagating large cracks.

This paper is organized in the following way. FPM formulations in Galerkin weak form are introduced in Section 2 and the interface debonding model is derived in Section 3. Then the algorithmic implementations of the proposed method are presented in Section 4. Finally, in section 5, the proposed method is validated through examples to predict fracture strength, and crack-growth trajectories of U-notched structures made of brittle materials. In Section 6, we complete this paper with some concluding remarks.

## 2 FPM formulations considering interface damage

For a linear elastic problem defined in the domain $\Omega$, the governing equations read as

$$\begin{cases} \sigma_{ij,j} + f_i = 0 \\ \varepsilon_{ij} = \frac{1}{2}(u_{i,j} + u_{j,i}) \quad \text{in } \Omega \\ \sigma_{ij} = D_{ijkl}\varepsilon_{kl} \end{cases} \quad (1)$$

where $\sigma_{ij}$, $\varepsilon_{ij}$, $u_i$, $f_i$, $D_{ijkl}$ denote stress field, strain field, displacement field, body force field, and the fourth-order stiffness tensor. Boundary conditions are given as,

$$\begin{cases} u_i = \overline{u}_i & \text{on } \Gamma_u \\ \sigma_{ij}n_j = \overline{t}_i & \text{on } \Gamma_t \end{cases} \quad (2)$$

where $n_j$ is the outward unit normal vector on the boundary portion $\Gamma_t$, $\overline{u}_i$ and $\overline{t}_i$ are the prescribed displacement and traction fields on the boundary portions $\Gamma_u$ and $\Gamma_t$, respectively. Since this work focuses on 2D problems, the variables take the value $i, j, k, l = 1, 2$.

In FPM, the problem domain is discretized by a series of points which are located inside the domain or on its boundary. Using these points, the problem domain is further partitioned into contiguous and non-overlapping subdomains. Details of discretization and partition methods have been introduced in detail in [31, 32].

In a subdomain $E_0$ for the point $P_0$, the trial function $\mathbf{u}^h = \left[u_1^h, u_2^h\right]^T$ is defined as

$$\mathbf{u}^h = \mathbf{u}^0 + \begin{bmatrix} \mathbf{h} & 0 \\ 0 & \mathbf{h} \end{bmatrix} \begin{bmatrix} \mathbf{a}_1 \\ \mathbf{a}_2 \end{bmatrix} \quad (3)$$

where $\mathbf{u}^0 = \left[u_1^0, u_2^0\right]^T$ is the displacement vector of $P_0$. For linear trial function,

$$\mathbf{h} = \begin{bmatrix} x_1 - x_1^0 & x_2 - x_2^0 \end{bmatrix}, \quad \mathbf{a}_i = \begin{bmatrix} \dfrac{\partial u_i}{\partial x_1} & \dfrac{\partial u_i}{\partial x_2} \end{bmatrix}^\mathrm{T}\bigg|_{P_0}, \quad i = 1, 2 \qquad (4)$$

where $x_i^0$ and $x_i$ are the coordinates of the point $P_0$ and an arbitrary point $P$ in $E_0$, respectively.

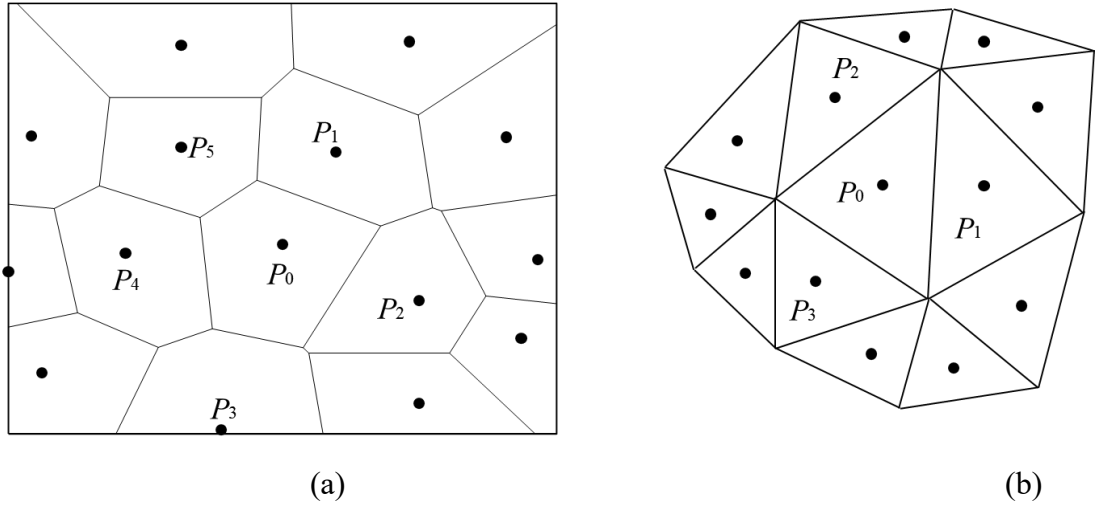

Fig. 1 Support of $P_0$ : (a) $P_0 \sim P_5$, (b) $P_0 \sim P_3$

In this work, the support of $P_0$ is defined as the points, $P_1 \sim P_m$, of neighboring subdomains, $E_1 \sim E_m$, which share the same boundary with the subdomain $E_0$ (see Fig. 1). Using the generalized finite difference method, the displacement gradient vector $\mathbf{a} = \begin{bmatrix} \mathbf{a}_1^\mathrm{T}, \mathbf{a}_2^\mathrm{T} \end{bmatrix}^\mathrm{T}$ in the subdomain $E_0$ can be derived as,

$$\mathbf{a} = \mathbf{C}\mathbf{u} \qquad (5)$$

where $\mathbf{C}$ serves as the gradient operator and $\mathbf{u} = [\mathbf{u}^{0\mathrm{T}}, \mathbf{u}^{1\mathrm{T}}, \ldots, \mathbf{u}^{m\mathrm{T}}]^\mathrm{T}$ collects displacement of all points in the support of $P_0$. Then, $\mathbf{u}^h$ can be expressed as

$$\mathbf{u}^h = \mathbf{N}\mathbf{u} \qquad (6)$$

where $\mathbf{N}$ is the shape function. The trial strain field $\boldsymbol{\varepsilon}^h$ and trial stress field $\boldsymbol{\sigma}^h$ can be expressed as

$$\boldsymbol{\varepsilon}^h = \mathbf{B}\mathbf{u}, \quad \boldsymbol{\sigma}^h = \mathbf{D}\boldsymbol{\varepsilon}^h = \mathbf{D}\mathbf{B}\mathbf{u} \qquad (7)$$

where $\mathbf{D}$ and $\mathbf{B}$ are the stiffness tensor and the strain shape function, respectively.

Detailed derivation and expression of the trial function can be found in [32]. Note that, though the generalized finite difference method is used in this work, the displacement gradient needed at the Fragile Points can also be derived using other alternative methods, such as differential quadrature[37].

In this work, FPM formulations are expressed in Galerkin weak form using the same type of trial function $\mathbf{u}^h$ and test function $\mathbf{v}$ in each subdomain as,

$$\sum_{E \in \Omega} \int_E \sigma_{ij}(\mathbf{u}) \varepsilon_{ij}(\mathbf{v}) d\Omega - \sum_{E \in \Omega} \int_{\partial E} \sigma_{ij}(\mathbf{u}) n_j v_i d\Gamma = \sum_{E \in \Omega} \int_E f_i v_i d\Omega \tag{8}$$

Since the trial and test functions are discontinuous, interior penalty numerical flux correction, which has been widely used in discontinuous Galerkin FEM, is employed in the Galerkin weak form of FPM to guarantee the consistency and stability of this method as discussed in [31]. Accordingly, the Galerkin weak form of FPM with numerical flux correction can be expressed as

$$\sum_{e \in \Omega} \int_E \sigma_{ij}(\mathbf{u}) \varepsilon_{ij}(\mathbf{v}) d\Omega - \sum_{e \in \Gamma_h} \int_e t_i^* [v_i] d\Gamma = \sum_{E \in \Omega} \int_E f_i v_i d\Omega + \sum_{e \in \Gamma_t} \int_e \bar{t}_i v_i d\Gamma \tag{9}$$

where $t_i^*$ is the interior penalty numerical flux defined on internal interface $\Gamma_h$ between each pair of neighboring subdomains.

**Remark 1:** The trial and test functions used in FPM are polynomial, thus the integrals of Galerkin weak form can be evaluated by means of simple Gauss Integrations of very low order.

In this work, the numerical flux is defined according to the IIPG method as,

$$t_i^* = \{\sigma_{ij}(\mathbf{u}) n_j\} - \alpha_{ij} [u_j], \quad \text{on } \Gamma_h \tag{10}$$

where $[\cdot]$ and $\{\cdot\}$ are the jump and average operators. For an arbitrary internal interface $e \in \partial E_1 \cap \partial E_2$ shared by subdomains $E_1$ and $E_2$, the operators act on arbitrary quantity $w$ as

$$[w] = w\big|_e^{E_1} - w\big|_e^{E_2}, \quad \{w\} = \frac{1}{2}\left(w\big|_e^{E_1} + w\big|_e^{E_2}\right) \tag{11}$$

In Eq. (10), $\alpha_{ij}$ is a second-order penalty tensor used to weakly enforce displacement continuity across interior interfaces. In this work, it is defined as a diagonal matrix

$$\alpha_{ii} = \lambda E / h_s \text{ and } \alpha_{ij} = 0, \quad \forall i \neq j \tag{12}$$

where $E$ is Young's modulus, $h_s$ is the characteristic length of subdomain, and $\lambda$ is an arbitrary penalty parameter. According to its definition, $\alpha_{ij}$ can also be considered as the interfacial stiffness, and thus the term $\alpha_{ij}[u_j]$ in Eq. (10) measures the contribution of separation on the interaction between neighboring subdomains. Therefore, the interior penalty numerical flux $t_i^*$ has the physical meaning as interaction between each pair of neighboring subdomains, or more specifically as the exact traction acting on each internal interface, considering both the stress in subdomains and the separation across internal interface.

**Remark 2:** In FPM, the trial functions are not independent between neighboring subdomains, and thus the DOFs of a subdomain are influenced by field solutions of neighboring subdomains. For this reason, using a small penalty parameter in the numerical flux is sufficient to weakly enforce continuity in FPM.

In Galerkin weak form of FPM Eq. (9), the first term on the left-hand side and the two terms on the right-hand side are similar to traditional Galerkin weak form of FEM, representing the strain energy and the external work done by prescribed body force and traction. The second term on the left-hand side of Eq. (9) is the interior penalty numerical flux correction and, considering the physical meaning of numerical flux $t_i^*$ discussed above, represents the energy stored at each internal interface. Note that, since the numerical flux is expressed in the IIPG form, a term has been dropped in Eq. (9) compared to the weak form FPM formulations in the previous study[31, 32], which was used to guarantee the symmetry of stiffness matrix.

Since the trial function used in FPM is discontinuous, interface damage and crack can be modeled explicitly within this method by means of following simple modifications.

Firstly, the undamaged interior penalty numerical flux $t_i^*$ should be replaced by damaged one $t_i^{*d}$. Then the Galerkin weak form of FPM considering interface damage can be expressed as

$$\sum_{e\in\Omega}\int_E \sigma_{ij}(\mathbf{u})\varepsilon_{ij}(\mathbf{v})d\Omega - \sum_{e\in\Gamma_h}\int_e t_i^{*d}[v_i]d\Gamma = \sum_{E\in\Omega}\int_E f_i v_i d\Omega + \sum_{e\in\Gamma_t}\int_e \bar{t}_i v_i d\Gamma \quad (13)$$

Damaged numerical flux $t_i^{*d}$ is determined based on the interface debonding model which will be introduced in the following section. Before damage initiation, damaged and undamaged numerical flux is identical $t_i^{*d} = t_i^*$, and the system remains elastic. When numerical flux of an interior interface reaches tensile strength, damage is initiated, and the interface starts to dissipate energy during damage evolution. In this case, the undamaged numerical flux $t_i^*$ is replaced by the damaged one $t_i^{*d}$ based on the interface debonding model. When the energy dissipated by a damaged interface reaches fracture energy, a crack (completed separation at the interface) is postulated to be formed and damaged numerical flux vanishes $t_i^{*d} = 0$.

Secondly, when an interior interface is damaged, the support of the pair of neighboring points should be modified to exclude each other from their support (see Fig. 2 corresponding to Fig. 1), which guarantees that the trial functions in this pair of neighboring subdomains become independent. In this way, the interaction between the neighboring subdomains is governed only by the damaged numerical flux $t_i^{*d}$ during the damage evolution, and the neighboring subdomains are completely separated when the crack is formed.

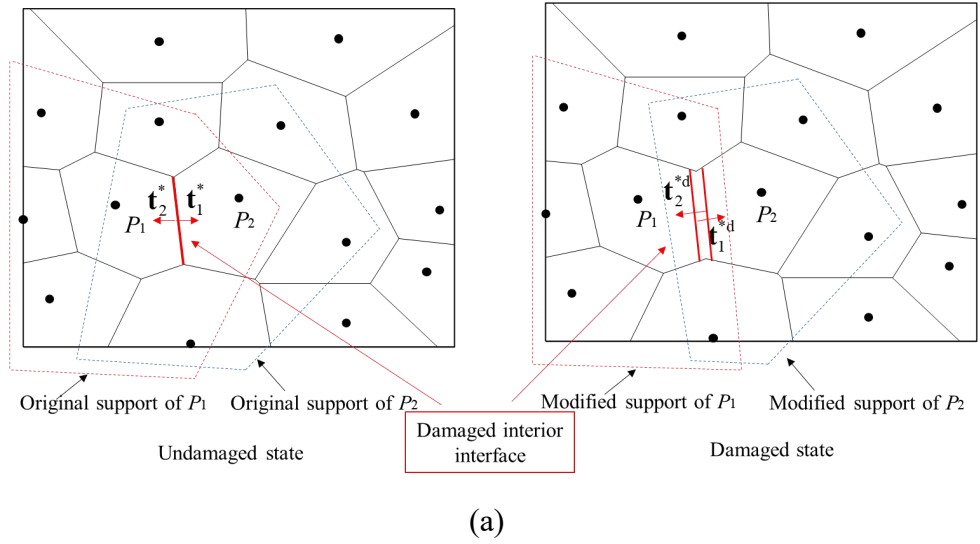

(a)

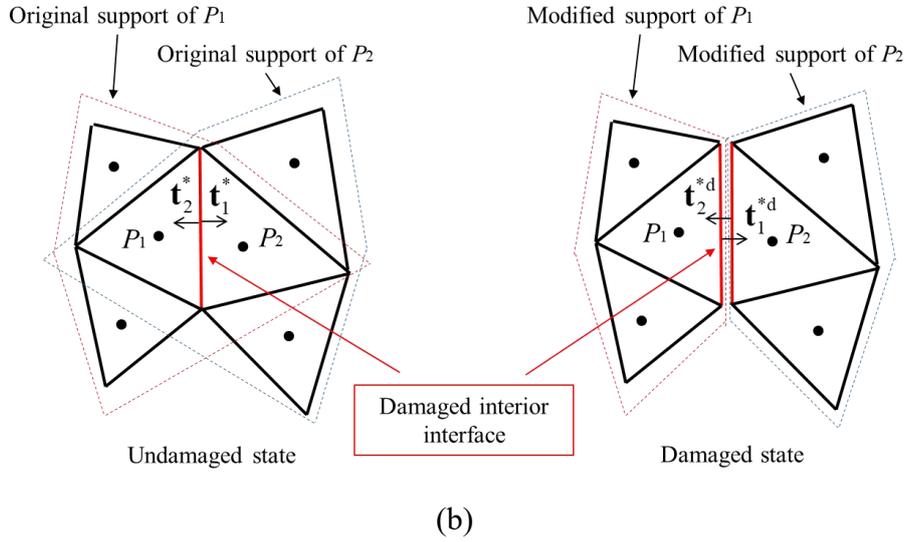

(b)

Fig. 2 Modification on support and numerical flux for interface damage

**Remark 3:** By correcting the numerical flux on an interior interface and the support of neighboring points (see Fig. 2), interface damage and crack can be simulated naturally and explicitly, with only simple and local modifications applied to FPM formulations. This is done without refining the mesh, inserting interface elements, and enriching shape functions, which are mostly required in numerical methods.

## 3 Formulations of the interface debonding model

In this work, an interface debonding model is derived based on a continuum damage framework, based on thermodynamics of irreversible processes and internal state variable theory [42, 43]. The mapping from the undamaged numerical flux $\mathbf{t}^*$ to the damaged one $\mathbf{t}^{*d}$ for an interior interface has been defined as,

$$\mathbf{t}^{*d} = (1-d)\mathbf{t}^* \tag{14}$$

where $d \in [0,1]$ is the interface damage scalar variable. The interfacial potential energy is defined as

$$\Psi(\boldsymbol{\delta}, d) = (1-d)\Psi^0(\boldsymbol{\delta}) \tag{15}$$

where $\Psi^0(\boldsymbol{\delta})$ is the equivalent interfacial elastic potential energy of undamaged interface reads

$$\Psi^0(\boldsymbol{\delta}) = \frac{1}{2}\boldsymbol{\delta}^{\mathrm{T}} \cdot \boldsymbol{\alpha} \cdot \boldsymbol{\delta} \tag{16}$$

Here, $\boldsymbol{\delta}$ is the thermodynamic conjugate vector of undamaged numerical flux $\mathbf{t}^*$ reads

$$\boldsymbol{\delta} = \boldsymbol{\alpha}^{-1} \cdot \mathbf{t}^* = \boldsymbol{\alpha}^{-1} \cdot \{\boldsymbol{\sigma} \cdot \mathbf{n}\} - [\mathbf{u}] \tag{17}$$

where $\boldsymbol{\alpha}$ is the interfacial second-order penalty tensor (or equivalent stiffness tensor) defined in Eq. (12).

Then the Clausius-Duhem inequality [44, 45] of interior interfaces takes the form

$$-\dot{\Psi} + \mathbf{t}^{*d} \cdot \dot{\boldsymbol{\delta}} \geq 0 \tag{18}$$

and the damaged numerical flux is derived as

$$\mathbf{t}^{*d} = \frac{\partial \Psi}{\partial \boldsymbol{\delta}} = (1-d)\frac{\partial \Psi^0}{\partial \boldsymbol{\delta}} \tag{19}$$

Substituting Eq. (15) into Eq. (18), one obtains the dissipation inequality as

$$\Psi^0 \dot{d} \geq 0 \tag{20}$$

which indicates that damage evolution is irreversible.

Defining the energy norm as

$$\tau(\boldsymbol{\delta}) = \sqrt{2\Psi^0(\boldsymbol{\delta})} \tag{21}$$

the damage initiation is characterized by means of the damage criterion based on $\boldsymbol{\delta}$ reads

$$f(\tau, r) = \tau - r \leq 0 \tag{22}$$

where $r$ is current damage threshold. Eq. (22) indicates that interface damage is initiated when the energy norm $\tau$ exceeds initial damage threshold $r_0$. Damage evolution is governed by

$$\dot{d} = \dot{\mu} h(\tau, d) \quad \text{and} \quad \dot{r} = \dot{\mu} \tag{23}$$

where $h$ controls damage evolution and $\dot{\mu}$ is an indicator for damage evolution according to Kuhn-Tucker conditions

$$\dot{\mu} \geq 0, \; f(\tau, r) \leq 0, \; \dot{\mu} f(\tau, r) = 0 \tag{24}$$

When $f < 0$, one obtains $\dot{\mu} = 0$ according to Eq. (24), and thus $\dot{d} = \dot{\mu} = \dot{r} = 0$, according to Eq. (23), which indicates that interface damage does not develop. When $\dot{\mu} > 0$, one obtains

$$f(\tau, r) = \dot{f}(\tau, r) = 0 \;\Rightarrow\; \dot{\mu} = \dot{r} = \dot{\tau} \tag{25}$$

and $\dot{d} > 0$ which indicates that damage evolution takes place.

Combining Eqs. (19), (21), (23) and (25), the following equation is derived,

$$\dot{\mathbf{t}}^{*d} = \mathbf{C}(\boldsymbol{\delta}) \cdot \dot{\boldsymbol{\delta}} \tag{26}$$

where $\mathbf{C}(\boldsymbol{\delta})$ is the tangent moduli expressed as

$$\mathbf{C}(\boldsymbol{\delta}) = (1-d)\frac{\partial^2 \Psi^0(\boldsymbol{\delta})}{\partial \boldsymbol{\delta}^2} - \frac{h(\tau, d)}{\tau}\frac{\partial \Psi^0(\boldsymbol{\delta})}{\partial \boldsymbol{\delta}} \otimes \frac{\partial \Psi^0(\boldsymbol{\delta})}{\partial \boldsymbol{\delta}} \tag{27}$$

In the currently presented interface debonding model, damage initiation and evolution are governed by the damage criterion Eq. (22) and evolution law Eq. (23), for which damage threshold $r$ and control function $h(\tau, d)$ have not been defined yet.

Fracture of brittle materials is commonly considered as being dominated by principal stress [14], so, substituting principal stress by normal traction, damage is initiated when the interface numerical flux reaches the interface strength $t^C$. Combining Eqs. (17) and (21), the initial damage threshold can be defined as

$$r_0 = \tau^C = \sqrt{2\Psi^0(\delta^C)} \tag{28}$$

where $\delta^C = (\lambda E / h_s)^{-1} t^C$.

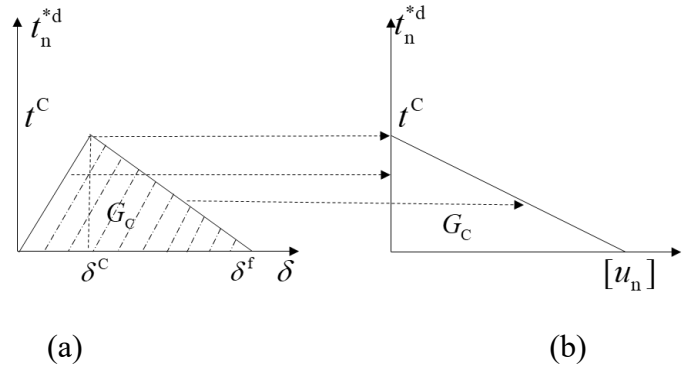

(a)            (b)

Fig. 3 Illustrative (a) $t_n^{*d} - \delta$ and (b) $t_n^{*d} - [u_n]$ curves for the cohesive interface debonding model

In conventional CZM, the constitutive behavior of interface is illustrated by traction-separation curve. Analogously, in the currently presented cohesive interface debonding model, the constitutive behavior of interface is illustrated as a curve for normal traction $t_n^{*d}$ versus its conjugate variable $\delta_n$, as plotted in Fig. 3a for linear softening case. According to Eq. (17), the equivalent curve for normal traction $t_n^{*d}$ versus normal separation $[u_n]$ is plotted in Fig. 3b. In this work, the softening law is defined to guarantee that the area enclosed by $t_n^{*d} - \delta_n$ curve and $t_n^{*d} - [u_n]$ curve is identical at any time step.

**Remark 4:** In the proposed interface debonding model, interior penalty numerical flux $\mathbf{t}^{*d}$, also named traction here, as well as its conjugate vector $\boldsymbol{\delta}$ are related to both separation across an interior interface and stress in neighboring subdomains. For this reason, though $\delta_n$ is larger than $0$, separation $[u_n]$ is closed to $0$ before damage is initiated, as shown in Fig. 3. This avoids using very large interfacial stiffness and the numerical problem commonly known as artificial compliance can be avoided naturally.

It is worthwhile mentioning that, since displacement continuity is enforced weakly in the proposed FPM approach, the interface separation $[u_n]$ is not exactly zero before damage initiation, and thus the numerical flux $t_i^{*d}$ is not strictly equal to the averaged interface traction $\{\sigma_{ij}(\mathbf{u})n_j\}$. Nevertheless, when the mesh is fine enough, the interface separation $[u_n]$ is very small. In this manuscript, the interface damage threshold $r$ is defined based on $t_n^{*d}$ reaching the interface strength, (see Eqs. (17), (21) and (22)), but it can also be defined based on the averaged interface traction $\{\sigma_{ij}(\mathbf{u})n_j\}$.

According to the currently proposed interface debonding model, the energy dissipated by a damaged interior interface can be calculated as

$$G = \int_0^{\delta^f} t_n^{*d} \mathrm{d}\delta \tag{29}$$

When the dissipated energy $G$ of an interface reaches the fracture energy $G_C$, the interface is completely damaged and a crack is initiated.

For a linear softening behavior as shown in Fig. 3a, the fracture energy can be expressed as

$$G_C = \frac{1}{2} t^C \delta^f \tag{30}$$

according to Eq. (29), and the following equation holds

$$\frac{\delta^C}{\delta} = 1 - (1 - \frac{\delta^C}{\delta^f})d \tag{31}$$

Substituting Eqs. (16), (21) and (30) into Eq. (31), the energy norm $\tau$ can be derived as

$$\tau = \frac{\tau^C}{1-Hd} \tag{32}$$

where $H = 1-(\tau^C)^2/2G_C$ is a constant related to critical energy norm $\tau^C$ and fracture energy $G_C$.

Combining Eqs. (23) and (32), the damage evolution control function $h$ can be derived as

$$h = \frac{(1-Hd)^2}{\tau^C H} \tag{33}$$

and the damage variable $d$ is explicitly determined as

$$d = \frac{1}{H}\left(1-\frac{\tau^C}{\tau}\right) \tag{34}$$

In summary, the proposed interface debonding model is developed based on a continuum damage framework, and the proposed method is thermodynamic consistent. When the type of softening law is determined, e.g. linear, only two parameters, interface strength $t^C$ and fracture energy $G_C$, are required to govern the initiation and evolution of interface damage. Note that, though the above derived interface debonding model shares some similarities with the one proposed by Versino et al. [45] which has been used in combine with DG-FEM, their model is derived based on the numerical flux $\mathbf{t}^{*d}$ and ours, instead, is derived based on the conjugate of numerical flux $\boldsymbol{\delta}$.

## 4 Algorithmic implementations

This section presents algorithmic implementations of the FPM with interface debonding model. Substituting Eqs. (6) and (7) into Eq. (13) and assembling over the discretized domain, Galerkin weak form of FPM can be expressed in the assembled from as

$$\mathbf{Ku} = \mathbf{F}_{ext} \tag{35}$$

where $\mathbf{K}$ is the assembled stiffness matrix, $\mathbf{u}$ is the assembly of all DOFs, and $\mathbf{F}_{ext}$ is the assembled external load vector. The assembled stiffness matrix $\mathbf{K}$ consists of the point stiffness matrix $\mathbf{K}_e$ and the boundary stiffness matrix $\mathbf{K}_h$, which can be expressed as

$$\mathbf{K}_e = \sum_{E \in \Omega} \int_E \mathbf{B}^T \mathbf{D} \mathbf{B} d\Omega \tag{36}$$

$$\mathbf{K}_h = \sum_{e \in \Gamma_h} -\int_e (1-d)\left([\mathbf{N}]^T \{\mathbf{n}_e \mathbf{D} \mathbf{B}\} - [\mathbf{N}]^T \boldsymbol{\alpha}[\mathbf{N}]\right) d\Gamma \tag{37}$$

where $\mathbf{n}_e$ is the matrix projecting stress into interface traction based on the outward normal vector.

The external load vector is assembled as

$$\mathbf{F}_{ext} = \sum_{E \in \Omega} \int_E \mathbf{N}^T \mathbf{f} d\Omega + \sum_{e \in \Gamma_t} \int_e \mathbf{N}^T \overline{\mathbf{t}} d\Gamma \tag{38}$$

where $\mathbf{f}$ and $\overline{\mathbf{t}}$ are the prescribed body forces and tractions.

The incremental form of assembled equilibrium formulations can be derived as

$$\hat{\mathbf{K}} \Delta \mathbf{u} = \mathbf{F}_{ext} - \mathbf{F}_{int} \tag{39}$$

in which the tangent stiffness matrix $\hat{\mathbf{K}}$ and internal load vector $\mathbf{F}_{int}$ reads

$$\hat{\mathbf{K}} = \sum_{E \in \Omega} \int_E \mathbf{B}^T \mathbf{D} \mathbf{B} d\Omega - \sum_{e \in \Gamma_h} \int_e \left([\mathbf{N}]^T \mathbf{C} \boldsymbol{\alpha}^{-1} \{\mathbf{n}_e \mathbf{D} \mathbf{B}\} - [\mathbf{N}]^T \mathbf{C}[\mathbf{N}]\right) d\Gamma \tag{40}$$

$$\mathbf{F}_{int} = \sum_{E \in \Omega} \int_E \mathbf{B}^T \boldsymbol{\sigma} d\Omega - \sum_{e \in \Gamma_h} \int_e [\mathbf{N}]^T \mathbf{t}^{*d} d\Gamma \tag{41}$$

Note that, any geometric nonlinearity has not been considered in this work, and thus all formulations have been constructed with respect to the material configuration.

Since the problem is nonlinear when interface damage is included, proper iteration method is required for implicit analysis. In this work, the Newton-Rapson (NR) method is used to solve the nonlinear problem. For the $n$th iteration of the time step $t + \Delta t$, Eq. (39) can be rewritten as

$$^{t+\Delta t}\hat{\mathbf{K}}^{(n)} \Delta \overline{\mathbf{u}}^{(n)} = {}^{t+\Delta t}\mathbf{F}_{ext}^{(n)} - {}^{t+\Delta t}\mathbf{F}_{int}^{(n)} \tag{42}$$

where $\Delta \overline{\mathbf{u}}^{(n)}$ denotes the displacement increment calculated in the $n$th iteration. Thus, the accumulated displacement increment of the time step $t + \Delta t$ is calculated as $^{t+\Delta t}\Delta \mathbf{u}^{(n)} = \sum_{i=1}^{n} \Delta \overline{\mathbf{u}}^{(i)}$ and the current displacement field is updated as

$$^{t+\Delta t}\mathbf{u}^{(n)} = {^t}\mathbf{u} + {^{t+\Delta t}}\Delta\mathbf{u}^{(n)}.$$

Detailed solution procedure of the NR iteration scheme is summarized in Algorithm 1.

**Algorithm 1**: FPM NR iteration scheme

---

0. Initiate variables

$$^{t+\Delta t}\mathbf{u}^{(0)} = {^t}\mathbf{u},\ ^{t+\Delta t}\boldsymbol{\sigma}^{(0)} = {^t}\boldsymbol{\sigma},\ ^{t+\Delta t}\boldsymbol{\varepsilon}^{(0)} = {^t}\boldsymbol{\varepsilon},\ ^{t+\Delta t}d^{(0)} = {^t}d$$

1. Calculate the displacement increment $\Delta\overline{\mathbf{u}}^{(n)}$ in the $n$th iteration

$$\Delta\overline{\mathbf{u}}^{(n)} = \left(^{t+\Delta t}\mathbf{K}^{(n)}\right)^{-1}\left(^{t+\Delta t}\mathbf{F}_{\text{ext}}^{(n)} - {^{t+\Delta t}}\mathbf{F}_{\text{int}}^{(n)}\right)$$

2. Update the accumulated displacement increment at the time step $t + \Delta t$

$$^{t+\Delta t}\Delta\mathbf{u}^{(n)} = {^{t+\Delta t}}\Delta\mathbf{u}^{(n-1)} + \Delta\overline{\mathbf{u}}^{(n)}$$

3. Update the trial displacement field at the time step $t + \Delta t$

$$^{t+\Delta t}\mathbf{u}^{(n)} = {^t}\mathbf{u} + {^{t+\Delta t}}\Delta\mathbf{u}^{(n)}$$

4. Calculate shape functions and update the trial strain field and trial stress field as

$$^{t+\Delta t}\Delta\boldsymbol{\varepsilon}^{(n)} = {^{t+\Delta t}}\mathbf{B}^{(n)}\ ^{t+\Delta t}\Delta\mathbf{u}^{(n)},\quad ^{t+\Delta t}\boldsymbol{\varepsilon}^{(n)} = {^t}\boldsymbol{\varepsilon} + {^{t+\Delta t}}\Delta\boldsymbol{\varepsilon}^{(n)}$$

$$^{t+\Delta t}\Delta\boldsymbol{\sigma}^{(n)} = {^{t+\Delta t}}\mathbf{D}^{(n)}\ ^{t+\Delta t}\Delta\boldsymbol{\varepsilon}^{(n)},\quad ^{t+\Delta t}\boldsymbol{\sigma}^{(n)} = {^t}\boldsymbol{\sigma} + {^{t+\Delta t}}\Delta\boldsymbol{\sigma}^{(n)}$$

5. Update interface-related quantities based on the interface damage update scheme which will be introduced in **Algorithm 2**.

6. Calculate new shape function if the support is modified, and update the stress and strain fields again similar to Step 4.

7. Check convergence using following criterion

$$\left\|\Delta\mathbf{u}^{(n)}\right\| / \left\|^{t+\Delta t}\Delta\mathbf{u}\right\| \leq tol$$

where $tol$ is the given tolerance and $\|\cdot\|$ is norm operator.

---

In this work, nonlinearity is related to interface damage, and the interface damage update scheme used in Step 5 of **Algorithm 1** is presented in the following.

**Algorithm 2**: Interface damage update scheme

---

0. Input interface parameters and known variables for calculation

   Interface parameters:

   interface tensile strength $t^\mathrm{C}$, fracture energy $G_\mathrm{C}$ and penalty tensor $\boldsymbol{\alpha}$

   Known variables;

   Trial stress field $^{t+\Delta t}\boldsymbol{\sigma}^{(n)}$, separation $^{t+\Delta t}[\mathbf{u}]^{(n)}$, damage variable $^{t+\Delta t}d^{(n)}$ and damage threshold $^{t+\Delta t}r$.

1. Calculate constants $\tau^\mathrm{C}$ and $H$, and the initial damage threshold $r_0$.

2. Calculate interface traction and its conjugate variables $^{t+\Delta t}\mathbf{t}^{*(n)}$ and $^{t+\Delta t}\boldsymbol{\delta}^{(n)}$ based on Eqs. (10), (14) and (17).

3. Calculate the trial energy norm $^{t+\Delta t}\tau^{(n)}$ according to Eq. (32)

4. Update the elastic-damage tangent moduli and the damage variable

   **if** $^{t+\Delta t}\tau^{(n)} \leq {}^{t+\Delta t}r$ **then**

   $$^{t+\Delta t}d^{(n)} = {}^{t+\Delta t}d^{(n-1)}$$

   $$^{t+\Delta t}\mathbf{C}^{(n)} = \left(1 - {}^{t+\Delta t}d^{(n)}\right)\left(\frac{\partial^2 \Psi^0}{\partial \boldsymbol{\delta}^2}\right)^{(n)}$$

   **else**

   $$^{t+\Delta t}d^{(n)} = \frac{1}{H}\left(1 - \frac{\tau^\mathrm{C}}{\tau^{(n)}}\right)$$

   $$^{t+\Delta t}\mathbf{C}^{(n)} = \left(1 - {}^{t+\Delta t}d^{(n)}\right)\left(\frac{\partial^2 \Psi^0}{\partial \boldsymbol{\delta}^2}\right)^{(n)} - \frac{h^{(n)}}{\tau^{(n)}}\left(\frac{\partial \Psi^0}{\partial \boldsymbol{\delta}} \otimes \frac{\partial \Psi^0}{\partial \boldsymbol{\delta}}\right)^{(n)}$$

   **end**

5. Update interface-related quantities and the support of points related to damaged interfaces as shown in Fig. 2.

## 5 Numerical examples

The proposed FPM formulations with interface debonding model has been implemented into in-house code. Some results are presented and discussed in this section. In Section 5.1, a mode I separation example has been employed to verify the proposed method for modeling cracks, by comparing FPM results with interface debonding model and FEM/CZM results. Then the proposed method has been applied to predict fracture strength (Section 5.2) and crack-growth trajectory (Section 5.3) for U-notched structures made of brittle materials, and the numerical results have been compared with experimental results available in literature as validations of the method.

5.1 Modeling of mode I separation

The mode I separation example has been defined as a rectangular plate under uniaxial loading as shown in Fig. 4. The plate was discretized by 8 points and partitioned into 8 subdomains. The loading was applied using displacement control. The material was assumed to be isotropic defined by the parameters Young's modulus $E=10$ and Poisson's ratio $v=0$. In the discretized model, interface damage was allowed only across the interior interfaces marked by red segments in Fig. 4. The tensile strength and fracture energy were specified as $t^C=1$ and $G_C=0.2$, respectively. In FPM, the penalty parameter $\lambda$, introduced in Eq. (12), should take the value within the range from $10^{-2}$ to $10^2$ as discussed in [32].

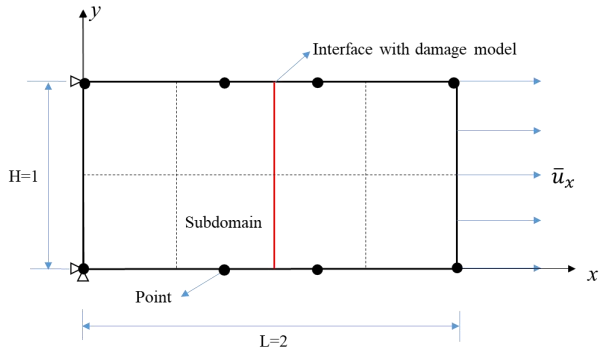

Fig. 4 FPM model for mode I separation

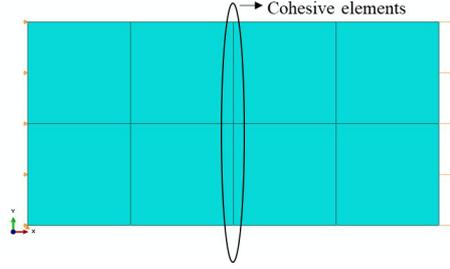

Fig. 5 FEM/CZM model for mode I separation

Using a similar setup, a FEM/CZM model with cohesive elements as shown in Fig. 5 has also been generated serving as a reference to verify FPM results. In the FEM/CZM model, before damage initiation, cohesive elements behave linearly [46], and large value of interfacial stiffness $K_{coh}$ should be taken to enforce a small separation. In this work, $K_{coh} = 10^4$ has been used in the FEM/CZM model illustrated in Fig. 5. As discussed in Remark 4, FPM allows the use of relatively small value of interfacial stiffness. In the FPM model shown in Fig. 4, the characteristic length was $h_s = 0.5$ and Young's modulus was $E = 10$, and the penalty parameter was taken as $\lambda = 1$, then the interfacial stiffness can be calculated as $\alpha = 20$ according to Eq. (12).

Taking the FEM results as reference, FPM results have been verified with respect to the deformation mode and stress field under various prescribed displacement, as summarized in Fig. 6, indicating that the proposed method can model mode I separation as accurately as FEM/CZM.

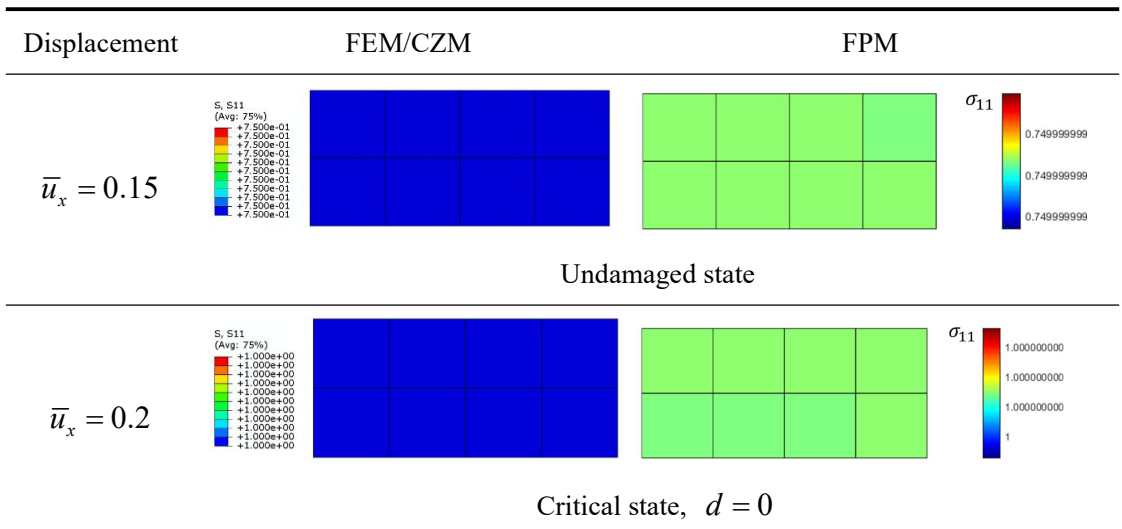

| | | |
|---|---|---|
| $\bar{u}_x = 0.25$ | | |

Damage state, $0 < d < 1$

| | | |
|---|---|---|
| $\bar{u}_x = 0.45$ | | |

Fractured state, $d = 1$

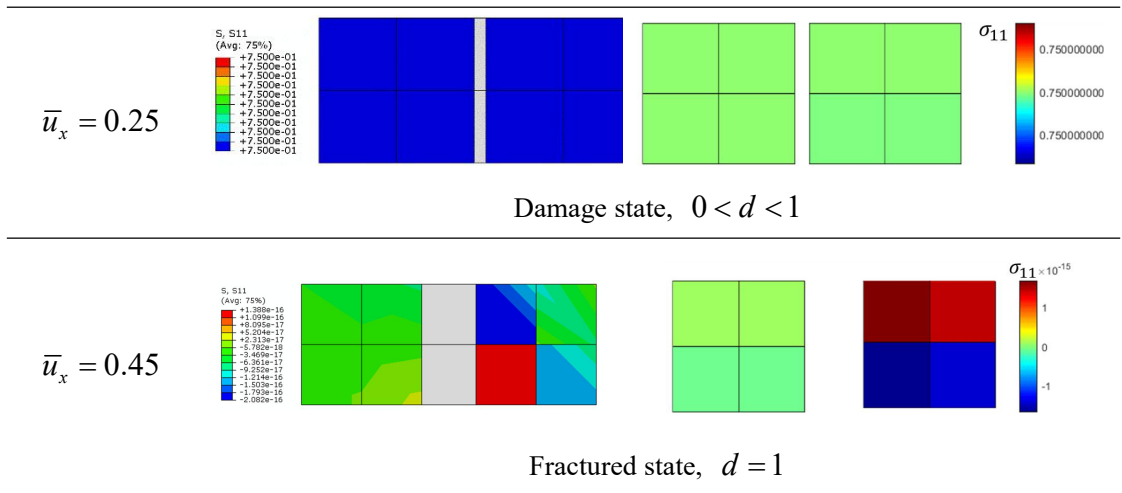

Fig. 6 Comparison of deformation mode and stress field obtained by FEM and FPM

In order to simulate mode II separation, the proposed FPM method has been applied to consider mixed-mode interface damage. As a verification, the uniaxial tensile loading in the FEM/CZM model in Fig. 4 and the FPM model in and Fig. 5 has been replaced by a pure shear loading. Similar to the mode I separation example, the comparison between FEM and FPM results given in Fig. 7 shows that, for mode II separation, the proposed FPM method with interface debonding model is also reliable.

| Displacement | FEM/CZM | FPM |
|---|---|---|
| $\bar{u}_y = 0.2$ | | |

Undamaged state

| | | |
|---|---|---|
| $\bar{u}_y = 0.4$ | | |

Critical state, $d = 0$

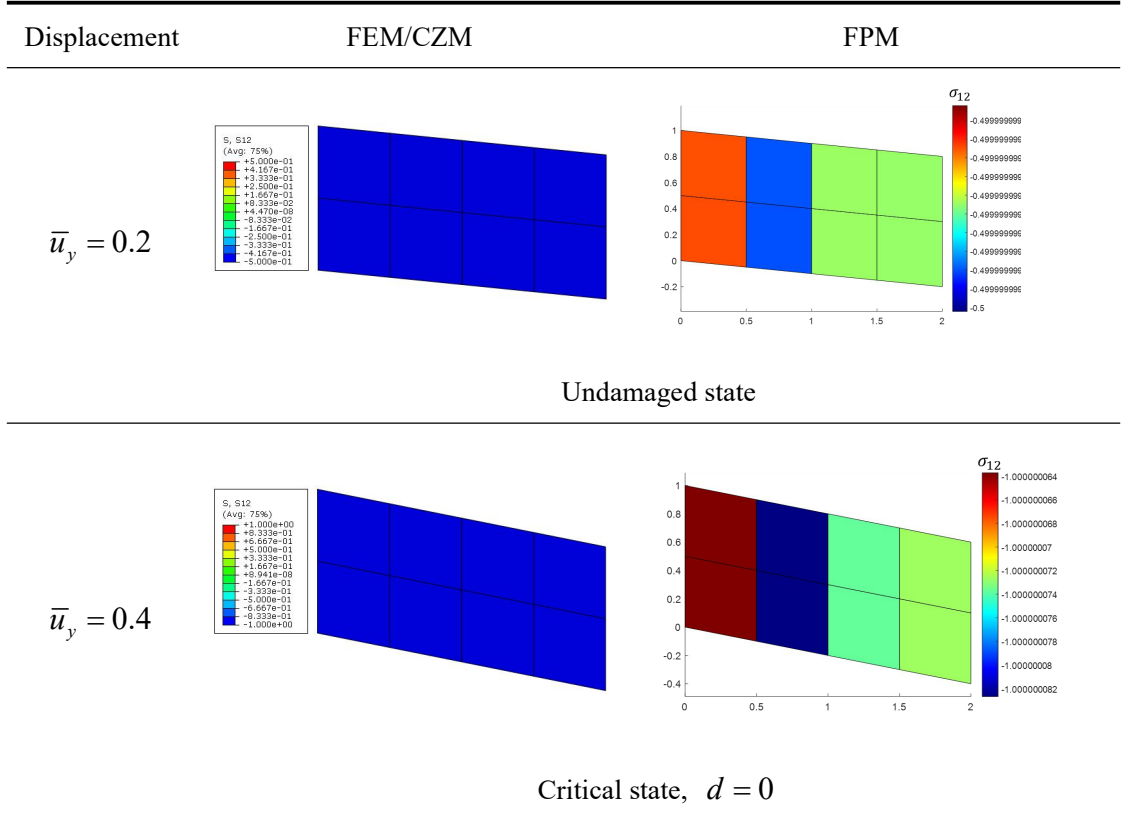

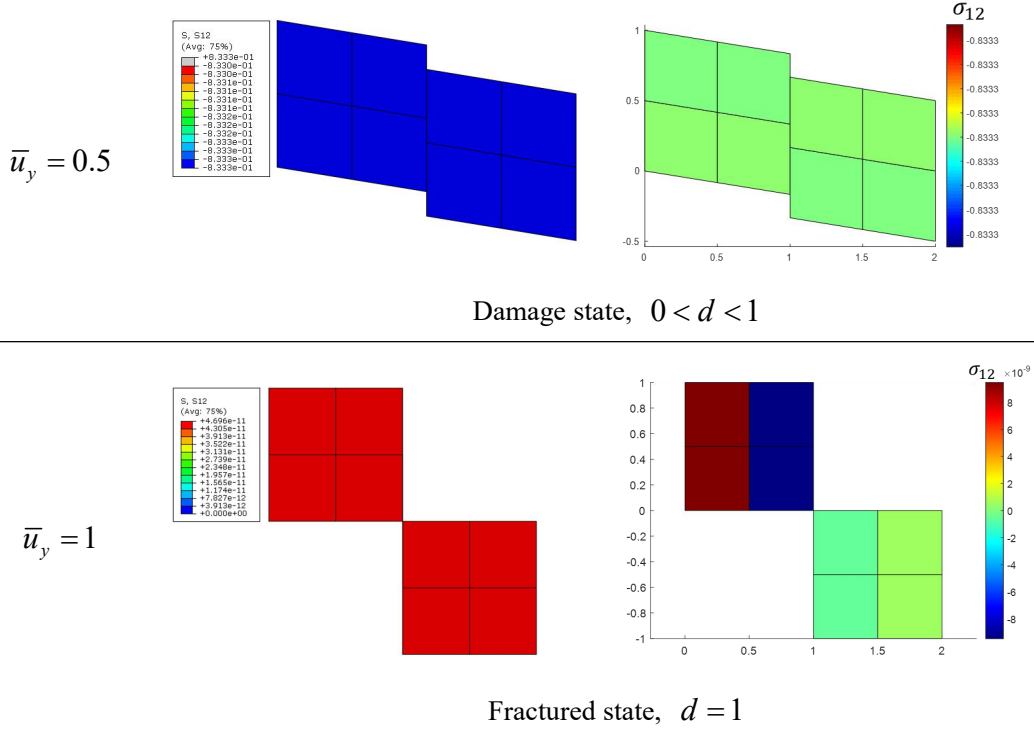

Fig. 7 Comparison of deformation mode and stress field obtained by FEM and FPM under a pure shear loading

As a next step, the influence of the interfacial stiffness, reads $K_{coh}$ in FEM and $\alpha$ in FPM has been studied by means of traction-separation curves plotted in Fig. 8a. It is observed that the shape of traction-separation curves for FEM is significantly influenced by the value of $K_{coh}$, and, in contrast, the curves for FPM are almost independent on the value of $\alpha$ when $\lambda$ is in the range from 1 to 100. The numerical results confirm that it is sufficient to use a small value of interfacial stiffness in FPM without losing the accuracy as discussed in Section 3. Additionally, Fig. 8b shows that, when the interface stiffness is not large enough, the FEM results are strongly mesh dependent even in the elastic regime, as the error increases when more cohesive elements along the loading direction are inserted in the FEM model. In contrast, the FPM results are not influenced by the number of interfaces as shown in Fig. 8b, verifying that the proposed FPM approach is less mesh dependent compared to the FEM approach coupled with CZM.

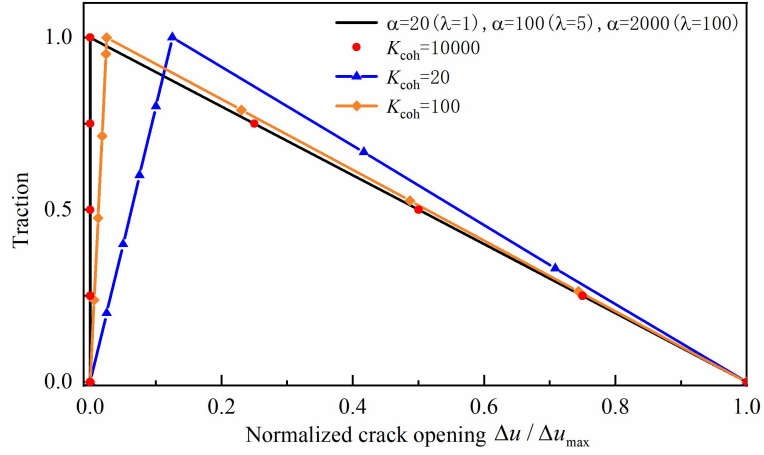

(a) traction-separation curves

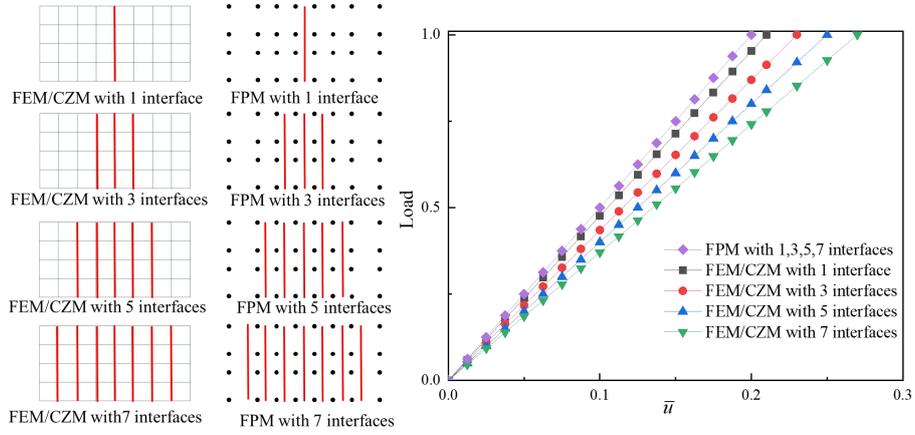

(b) load-displacement curves

Fig. 8 Numerical results illustrating the influence of (a) interface stiffness and (b) number of interfaces for FEM/CZM and FPM models

It is worth mentioning that the eight-point FPM model in Fig. 4 has 16 DOFs, and the eight-element FEM model in Fig. 5 has 18 nodes and 36 DOFs. For this reason, one can conclude that the proposed FPM with interface debonding model requires less DOFs than FEM with cohesive elements, and thus it can be expected as an efficient alternative method for modeling mode I fracture problems.

In this work, the computations were performed in a Windows desktop machine with an Intel(R) Core(TM) i7-6700 quad-core 3.4GHz processor. FEM results were obtained from commercial software package Abaqus/Standard and FPM calculations were carried out by self-written codes in MATLAB. For the mode I separation example (see Fig. 6), the CPU time of FEM and FPM were 0.5s and 1.25s, respectively. Note that, the self-written FPM codes are written in MATLAB, which are expected to have higher efficiency after standard code development using C++ or fortran.

5.2 Predicting fracture strength of U-notched specimens

In this section, the proposed FPM with interface debonding model has been applied to predict the fracture strength of U-notched structures under mode I loading with various root radii and brittle materials.

For each U-notched structure, the peak load, the structure can bear, was predicted by FPM simulation using the proposed interface debonding model, and thus the nominal stress $\sigma_{NC}$ at peak load can be obtained from the FPM result. In order to compare with experimental results from literature, the critical values of the generalized stress intensity factor was evaluated according to [4] reads

$$K_C^U = K_t \sigma_{NC} \sqrt{\pi \frac{R}{4}} \tag{43}$$

where $K_t$ is the notch stress concentration factor and $R$ is the U-notch radius. In this work, $K_t$ was obtained by FPM simulations without damage. Note that the critical generalized stress intensity factor $K_C^U$ is normally used as the measure of fracture toughness of U-notched structures in structural safety assessment.

The material parameters required by the proposed FPM approaches with interface debonding model are Young's modulus $E$, Poisson's ratio $\nu$, tensile strength $f_t$ and fracture energy $G_C$. In this work, it was assumed that the interfacial tensile strength $t^C$ equaled to the material's tensile strength $f_t$ Within the realm of LEFM [4, 11], the fracture energy $G_C$ is evaluated as

$$G_C = \frac{K_{IC}^2}{E/(1-\nu^2)} \tag{44}$$

where $K_{IC}$ is the material fracture toughness. For the brittle materials considered in this work, material parameters have been summarized in Table 1.

Table 1 Summary of material parameters used in this work

| Material | $E$ (GPa) | $v$ | $f_t$ (MPa) | $K_{IC}$ (MPa·m$^{1/2}$) | $G_C$ (J/m$^2$) |
| --- | --- | --- | --- | --- | --- |
| PMMA −60°C [11] | 5.05 | 0.4 | 128 | 1.7 | 480 |
| Al$_2$O$_3$ [1, 4, 47] | 373 | 0.2 | 297 | 3.8 | 37 |
| Y-PSZ [48] | 207 | 0.2 | 425 | 5.9 | 161 |

Double-edge U-notched tensile (DEUNT) specimens shown in Fig. 9a have been used to validate the proposed methods for predicting fracture strength. To reduce computational cost, the simplified model (see Fig. 9b) has been used in FPM simulations based on the Saint-Venant principle and symmetry of the problem. In this work, a convergence study has been carried out to determine the optimal size of simplified model as $B/R = 5$ and $L/R = 10$. The mesh size of FPM model shown in Fig. 9b was also determined based on a convergence study. Note that, in this work, the geometric models were discretized using FEM software ABAQUS, and then the FEM elements were converted into FPM subdomains. The effectiveness of the size and mesh of simplified model has been verified, by means of the specimen with PMMA (−60°C) and the root radius $R = 2$, as using full size and simplified DEUNT models resulted in similar stress concentration factor $K_t$.

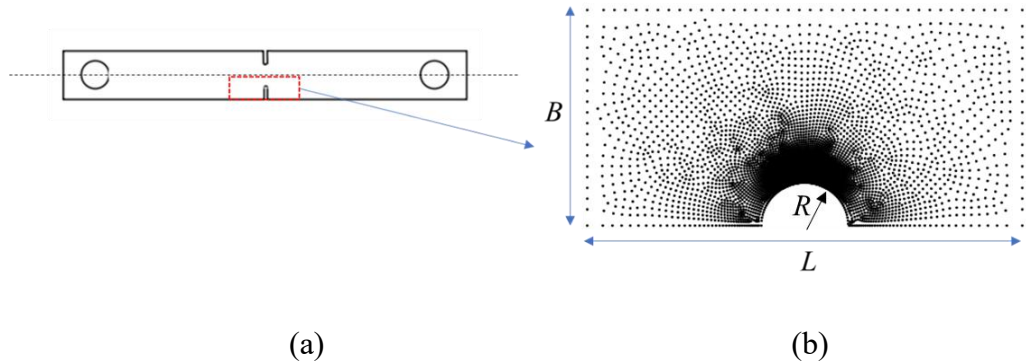

(a)  (b)

Fig. 9 Illustrative (a) geometry and (b) simplified FPM model of DEUNT specimens

For various brittle materials and root radii, the U-notch fracture toughness $K_C^U$ predicted by FPM with interface debonding model has been compared with experimental results available in literature [11, 47, 48]. The good agreement between numerical and experimental results summarized Fig. 10 indicates that the proposed method can predict fracture toughness as well as fracture strength of U-notched structures accurately regardless of the chosen materials and geometries.

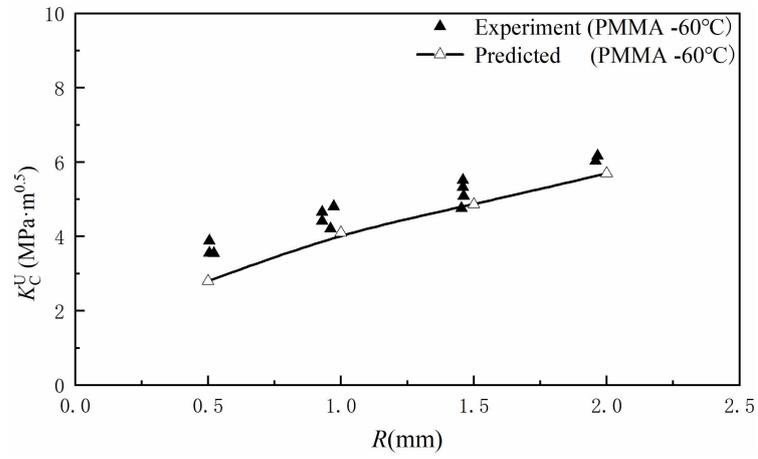

(a)

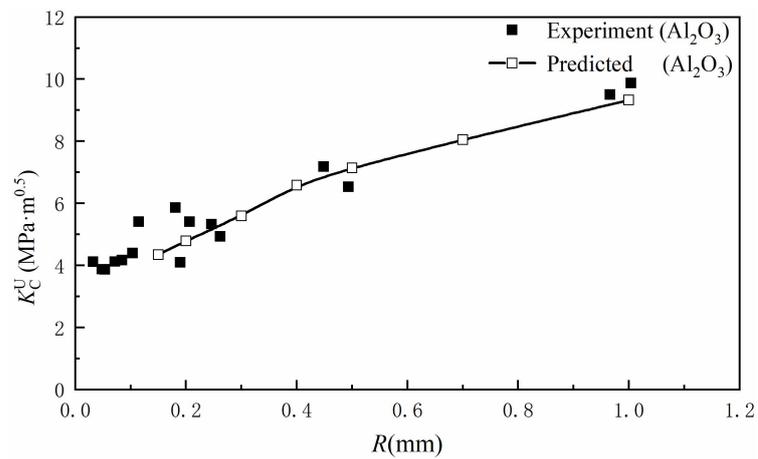

(b)

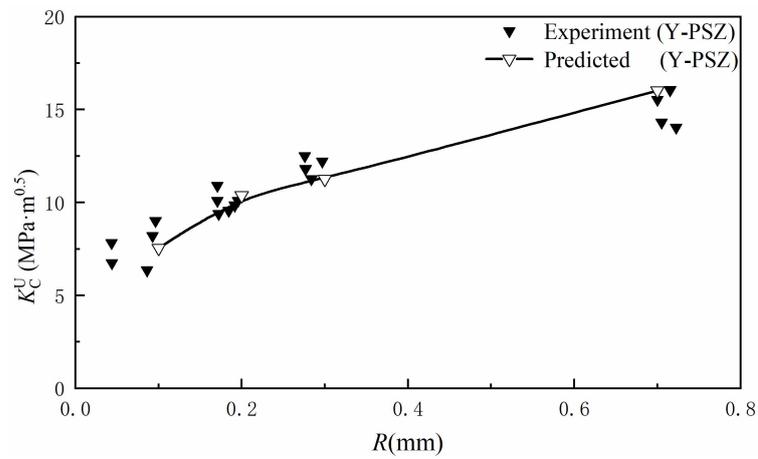

(c)

Fig. 10 Comparison on experimentally- and numerically-obtained U-notch fracture

toughness $K_C^U$ for (a) PMMA -60°C [11], (b) Al$_2$O$_3$ [47], and (c) Y-PSZ [48]

5.3 Predicting crack initiation from U-notches and the peak load with different length of grown crack

Predicting crack-growth trajectory is demanded during the design and application of engineering structures, since it helps to explore fracture mechanism and optimize structural design. In the proposed FPM approach, interface damage is characterized across each interior interface using the interface debonding model, presented in Section 3, in which the damage behavior is defined mainly based on tensile strength and fracture energy of the materials. Consequently, it is convenient and efficient to locate the crack initiation at the critical interface and track the crack propagation along the optimal path with minimum energy dissipated by crack surfaces.

To demonstrate the proposed method's capability on predicting crack-growth trajectory, two examples have been defined as shown in Fig. 11. The first example was a DEUNT specimen with PMMA (-60°C) and $R=2$, for which the experimental results from [11] were taken as references. The second example is a U-notched Brazilian disc (UNBD) specimen made of graphite under compression load. Two angles of U-notches, $0°$ and $10°$ have been considered in the simulations and the experimental results from [21, 49] were used for validation. In the simulations, the material parameters for PMMA (-60°C) were taken from Table 1 and those for graphite were prescribed as $E=8.05\text{GPa}$, $v=0.2$, $f_t=27.5\text{MPa}$ and $G_C=119\text{J/m}^2$ [49]. On the potential crack path, the location of points and partition of the problem domain are shown in Fig. 11c.

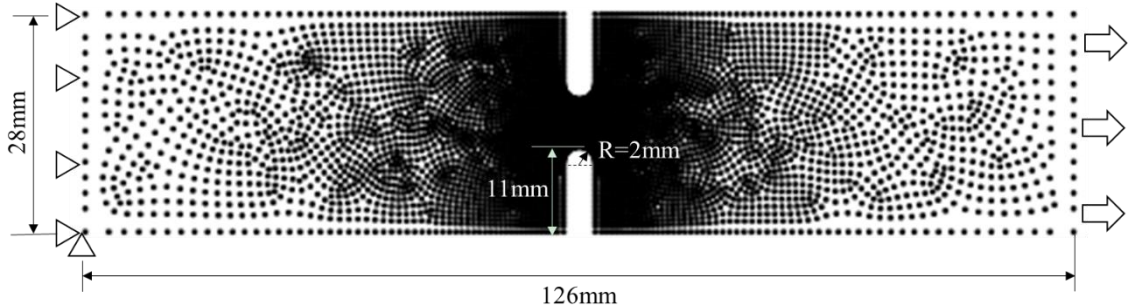

(a)

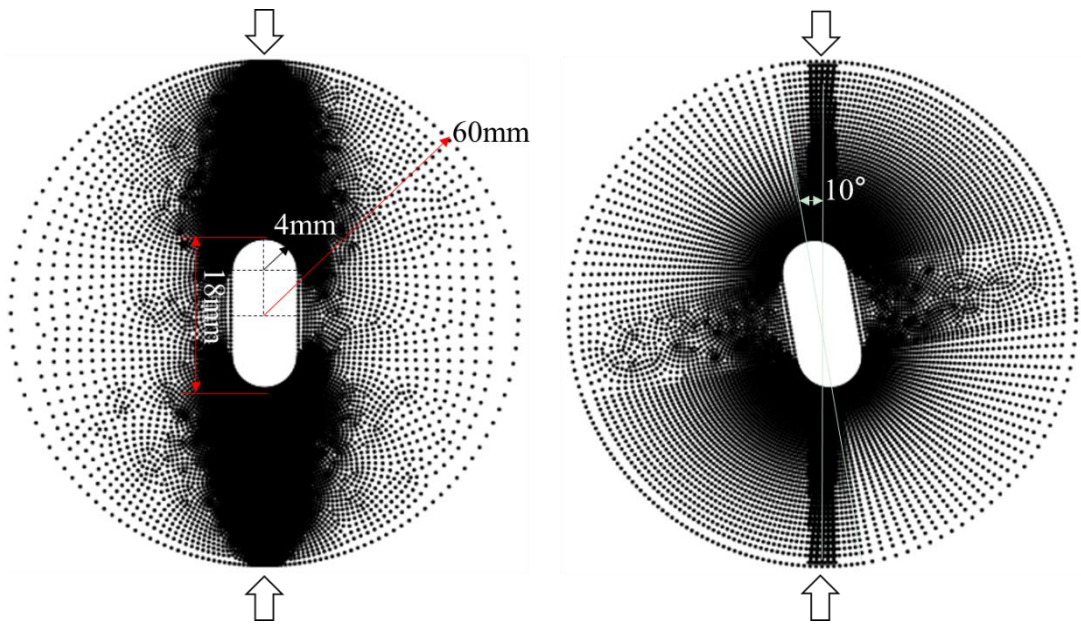

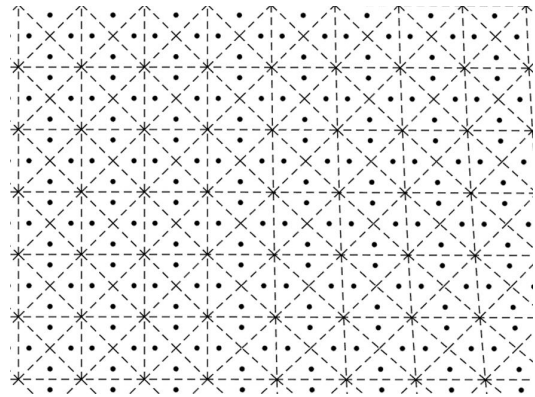

(c)

Fig. 11 FPM models of (a) the DEUNT, (b) UNBD specimens as illustrative examples for fracture trajectory prediction and (c) the location of Fragile Points, and the partition of the problem on the potential crack path

Modeling structural failure of notched specimens made of brittle materials has been a well-known challenging task for numerical methods, since, after reaching the peak load, the crack propagation is normally an unstable process [1, 2]. Thus it is difficult for an implicit method to converge and some additional techniques are needed, such as the arc-length method [50].

In this work, an implicit scheme has been developed to model the crack initiation and propagation process of U-notched specimens made of brittle materials.

For a first step, the peak load of a structure has been determined by means of the incremental method when the critical interior interface $\Gamma_{crit}$ reached the critical state. Here, $\Gamma_{crit}$ refers to the interface with the largest damage variable (less than 1) in the FPM model at the time step $t$. The critical state is defined as, when any further small load increment (smaller than the prescribed tolerance) is applied, the damage variable of $\Gamma_{crit}$ reaches 1 during iteration and the NR scheme fails to converge.

Taken the DEUNT model as an example, Fig. 12 shows the interface damage states at various loading steps. It can be observed that more interfaces are damaged under higher load, and the crack initiated across the interface with maximum damage variable (marked as Max in Fig. 12) under the maximum prescribed displacement, reads 0.5mm for the DEUNT model.

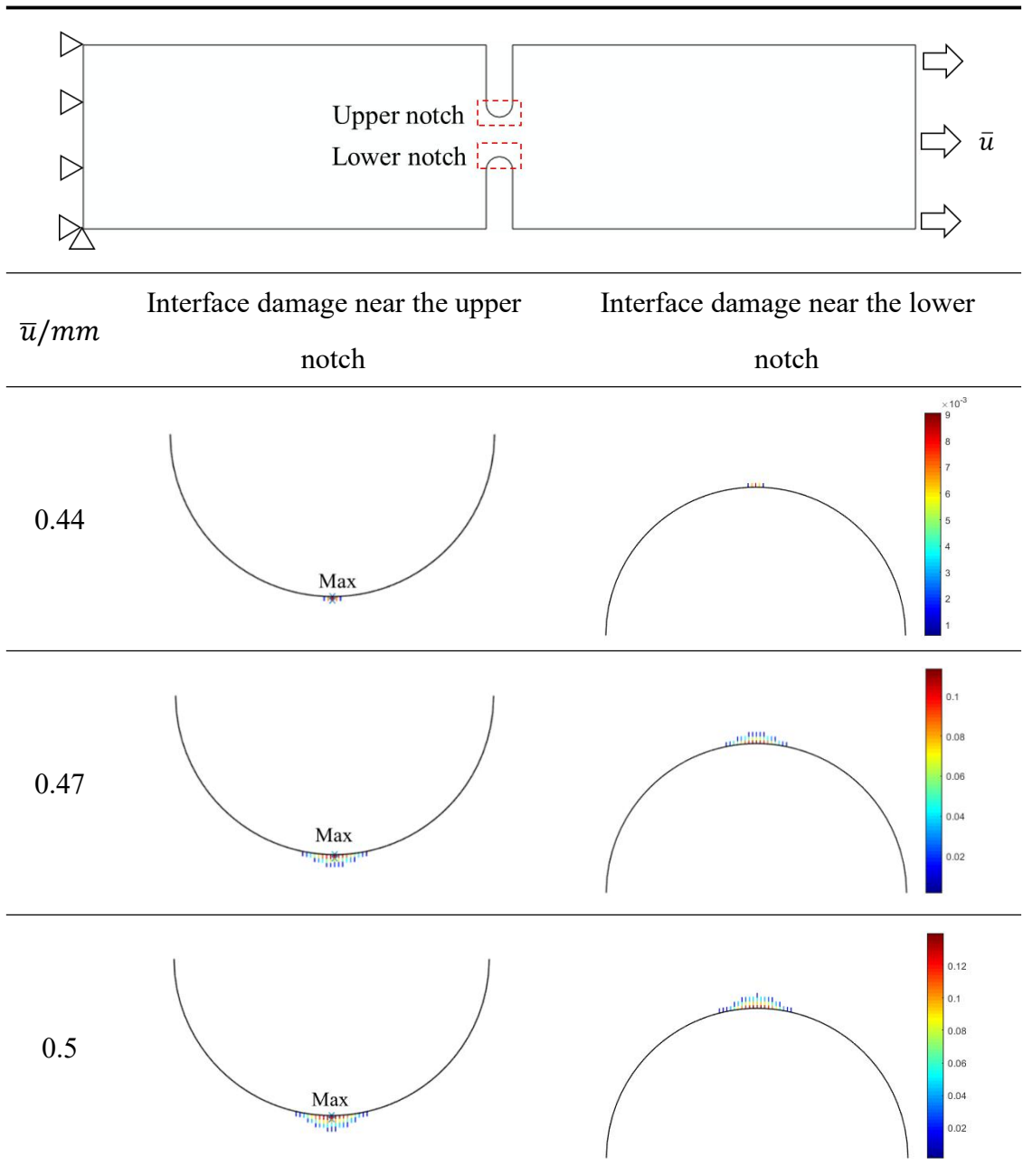

Fig. 12 Interface damage state of the DEUNT model

For a second step, a crack was inserted across the critical interface $\Gamma_{\text{crit}}$ and the damage variable of other interior interfaces in the FPM model remained unchanged. Then, by applying the external load incrementally again, the peak load and the crack development of the cracked structure was determined following the first step. The crack-growth trajectory can be predicted by repeating above procedure until the structure failed.

Fig. 13 shows the peak load-displacement curve of the DEUNT example which was determined in the following way. Firstly, the peak load $P = 6800\text{N}$ was determined for the intact specimen marked as red point in Fig. 13. Introducing an initial crack with a length of $0.05\text{mm}$, the peak load $P = 3763\text{N}$ was determined for the cracked specimen marked as blue point in Fig. 13. It is worthwhile to mention that the intact and cracked DEUNT specimens reached peak load at the prescribed displacement $0.5\text{mm}$ and $0.31\text{mm}$. The results indicated that equilibrium cannot be achieved after the first crack was initiated under the prescribed displacement $0.5\text{mm}$ which proved the instability of the consequent crack propagation. Using the implicit scheme introduced above, a series of peak loads with respect to crack length was obtained illustrated as the points in Fig. 13 and the peak load-displacement curve was plotted by connecting these points. In this way, one can predict the peak load as well as the crack development of U-notched structures with a crack grown from the root with various length.

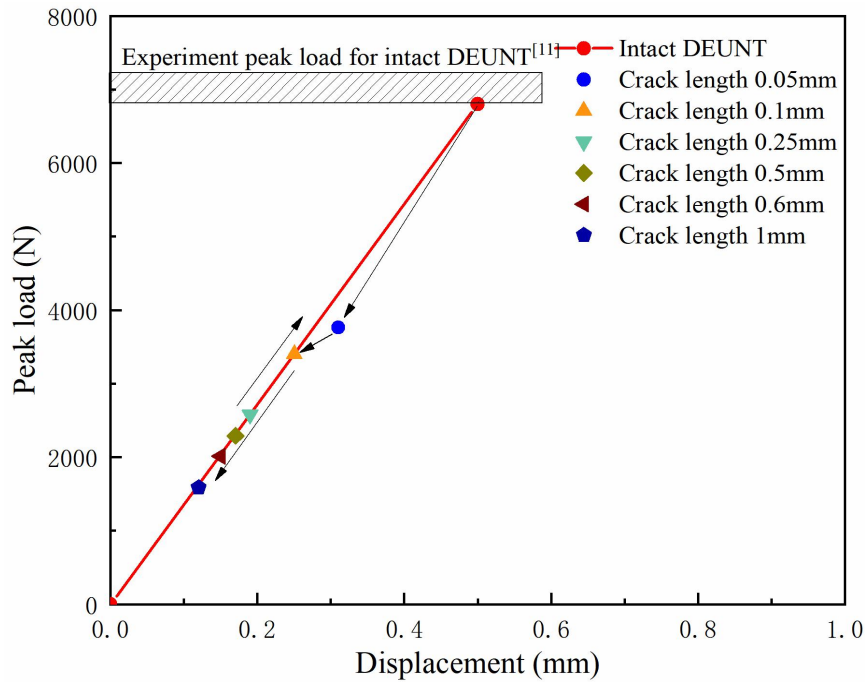

Fig. 13 Peak load-displacement curve of the DEUNT example

Fig. 14 illustrated the crack growth of DEUNT and UNBD specimens obtained from FPM simulations. For the DEUNT specimen, the crack propagated along the symmetry plane of the notch (see Fig. 14a), which conformed to the expected result. For the two UNBD specimens, the crack trajectory predicted by FPM agreed well with the experimental results available in [21]. The comparison in Fig. 14 validates the proposed incremental scheme based on FPM approach for predicting fracture trajectory of U-notched specimens made of brittle materials. The numerical results agree well with experimental observations for specimens under either mode I (see Figs. 14(a) and (b)) or mixed mode loading (see Fig. 14(c)).

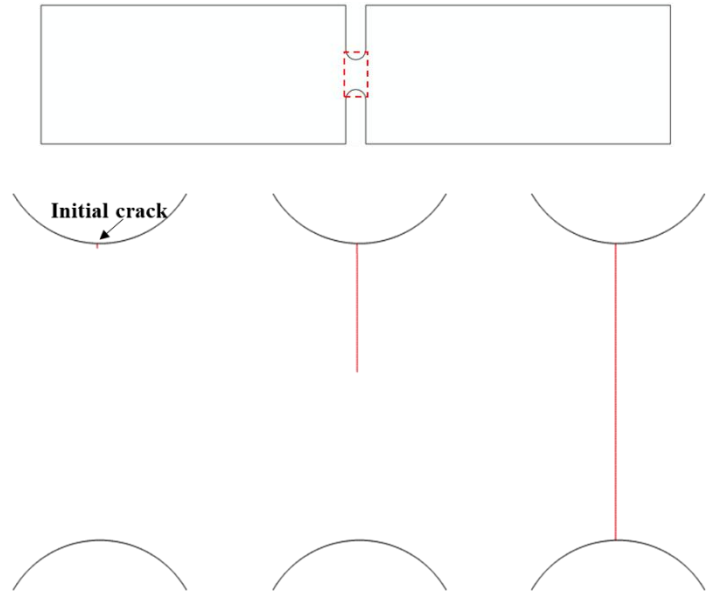

(a) DEUNT

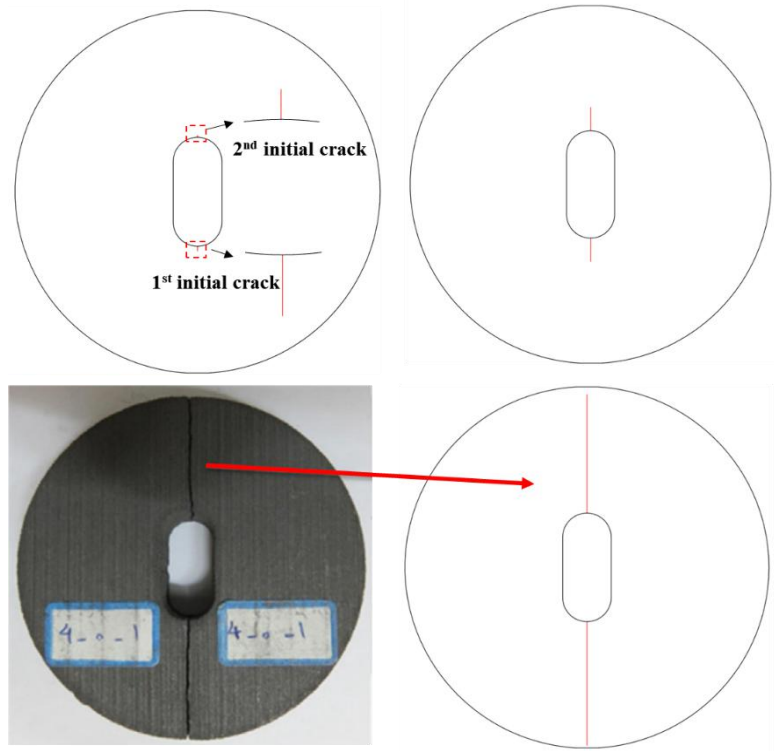

(b) UNBD (load angle=0°)

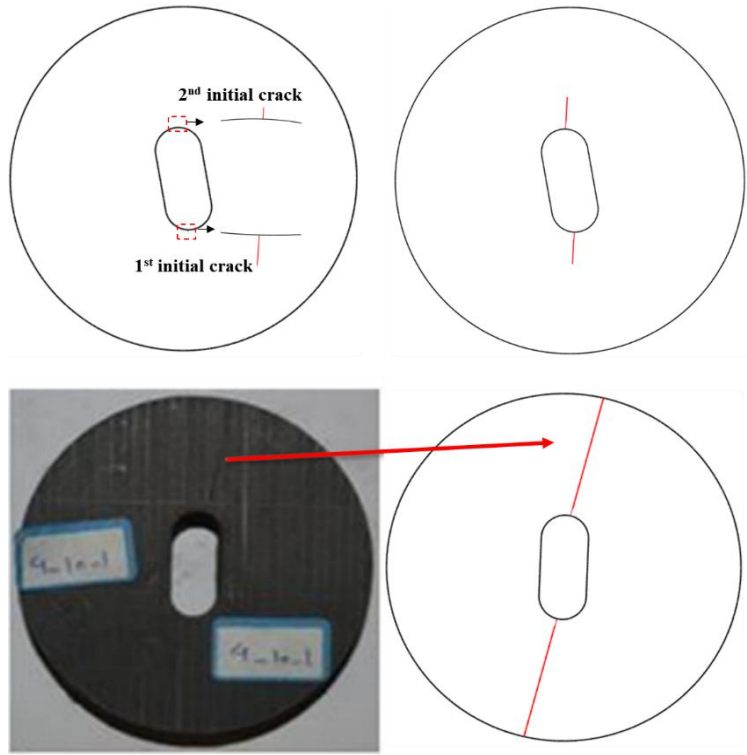

(c) UNBD (load angle=10°)

Fig. 14 Crack-growth trajectory prediction for DEUNT and UNBD specimens using FPM and the validation using experimental results [21]

One may argue that the crack paths shown in Fig. 14 were however predicted using structured mesh (see Fig. 11c). For a next step, the DEUNT and UNBD models were discretized by unstructured meshes (see Fig. 15). By comparing numerical results in Fig. 14 and Fig. 15, one can conclude that similar crack paths have been predicted by FPM using either structured or unstructured meshes, which indicates that the prediction of crack growth using FPM are mesh independent.

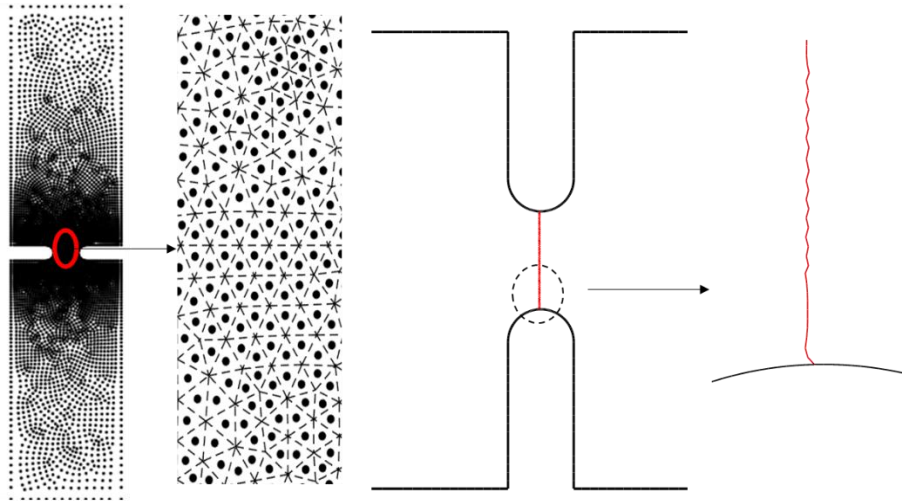

(a) DEUNT

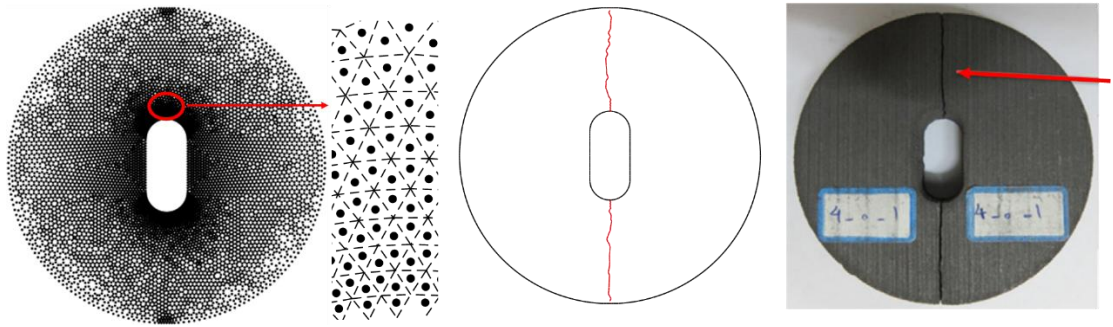

(b) UNBD (load angle=0°)

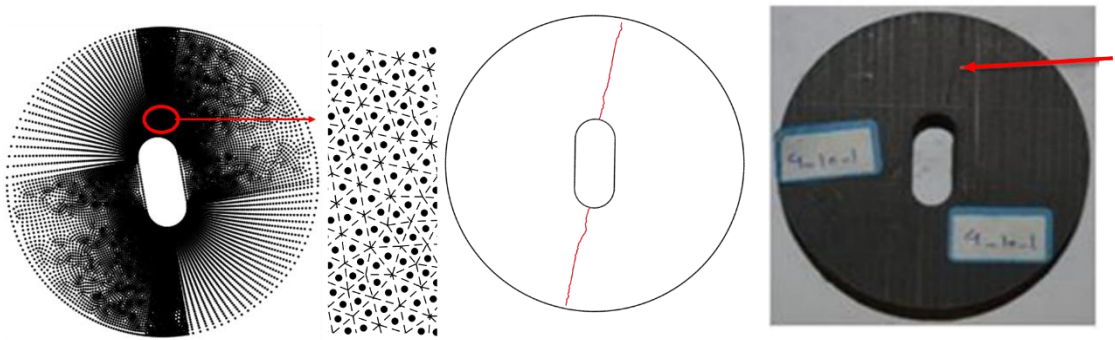

(c) UNBD (load angle=10°)

Fig. 15 Fracture trajectory prediction for DEUNT and UNBD specimens using unstructured mesh

Compared to existing numerical methods, the proposed FPM approach allows the use of consistent criteria and schemes for predicting fracture strength and crack-growth trajectory of U-notched structures made of brittle materials, and regenerating meshes or inserting additional DOFs are not required during the simulation. Thus, this method can be expected to serve as a powerful tool for predicting fracture behavior thanks to its effectiveness and convenience on modeling cracks.

**6 Conclusion**

In this study, basic formulations and algorithmic framework of the FPM approach with an interface debonding model have been proposed for damage assessment, such as predicting fracture strength and crack growth process of U-notched structures. FPM is a meshless method based on Galerkin weak form, and the trial functions are point-based, polynomial and discontinuous. To guarantee consistency of the method, interior penalty numerical flux corrections have been employed in the weak form formulations. In this work, the numerical flux has been introduced based on the IIPG method, and it has an explicit physical meaning as the exact traction acting on each interior interface. Then the damaged numerical flux is defined and introduced into Galerkin weak form. In this way, the crack initiation and propagation can be modeled explicitly by adjusting the support of points and updating the damaged numerical flux, without the need of remeshing or deleting elements.

Furthermore, an interface debonding model has been derived in a continuum damage framework to characterize the interface damage behavior. In this model, the damage variable has been defined mapping the undamaged numerical flux into the damaged one. The evaluation of damage variable is governed by the damage initiation criterion and softening law (assumed to be linear in this work), which are related to tensile strength and fracture energy, respectively. The interface damage is initiated when the numerical flux on the interface reaches tensile strength, and the crack is formed when the dissipated energy by the interface equaled to fracture energy. Combining the developed interface debonding model with FPM, since the point-based trial functions in neighboring subdomains are not independent, it is not necessarily to use large interfacial stiffness and thus the numerical problem commonly known as artificial compliance can be avoided naturally.

To assess the capability of the proposed method for solving practical engineering problems, FPM with the presently proposed interface debonding model has been applied to predict fracture strength and crack grown from U-notches in structures made of brittle materials. Good agreement between numerical and experimental results have been observed for the illustrative examples. Thus, this method has great potential for predicting damage and fracture behavior of notched structures in engineering practice.

In our future work, the proposed FPM method will be extended to study complex crack initiation and growth in the microstructure of composite material. In order to capture the nonlinear material response before crack initiation, more comprehensive constitutive model should be implemented to govern the behavior of subdomains, and a thorough study is needed to explore the coupling effect between the nonlinear constitutive model of subdomains and the debonding model of interior interfaces. Additionally, the proposed FPM formulations will be extended for solving dynamic problems, and, due to its strength in modeling cracks, this method will be applied to predict complex dynamic fracture behavior of brittle materials, such as crack branching in ceramics.

**Acknowledgements**

The first four authors acknowledge the support of the National Natural Science Foundation of China (12072011 and 12102023).

**References**

1. Gómez FJ, Elices M. A fracture criterion for blunted V-notched samples. *International Journal of Fracture*. 2004;127(3):239-64. doi: 10.1023/B:FRAC.0000036832.29429.21.
2. Gómez FJ, Elices M, Valiente A. Cracking in PMMA containing U - shaped notches. *Fatigue & Fracture of Engineering Materials and Structures*. 2000;23(9). doi: 10.1046/j.1460-2695.2000.00264.x.
3. Lazzarin P, Berto F. Some Expressions for the Strain Energy in a Finite Volume Surrounding the Root of Blunt V-notches. *International Journal of Fracture*. 2005;135:161-85. doi: 10.1007/S10704-005-3943-6.
4. Gómez FJ, Guinea GV, Elices M. Failure criteria for linear elastic materials with U-notches. *International Journal of Fracture*. 2006;141(1):99-113. doi: 10.1007/s10704-006-0066-7.
5. Ayatollahi M, Torabi A. Brittle fracture in rounded-tip V-shaped notches. *Materials & Design*.

2010;31(1):60-7. doi: 10.1016/J.MATDES.2009.07.017.

6. Ayatollahi M, Torabi A. Tensile fracture in notched polycrystalline graphite specimens. *Carbon*. 2010;48(8):2255-65. doi: 10.1016/J.CARBON.2010.02.041.

7. Torabi A. Fracture assessment of U-notched graphite plates under tension. *International Journal of Fracture*. 2013;181(2):285-92. doi: 10.1007/S10704-012-9799-7.

8. Lazzarin P, Berto F, Radaj D. Fatigue‐relevant stress field parameters of welded lap joints: pointed slit tip compared with keyhole notch. *Fatigue & Fracture of Engineering Materials & Structures*. 2009;32(9):713-35. doi: 10.1111/j.1460-2695.2009.01379.x.

9. Lazzarin P, Lassen T, Livieri P. A notch stress intensity approach applied to fatigue life predictions of welded joints with different local toe geometry. *Fatigue & Fracture of Engineering Materials & Structures*. 2003;26(1):49-58. doi: 10.1046/j.1460-2695.2003.00586.x.

10. Yosibash Z. Failure criteria for brittle elastic materials In: Singularities in Elliptic Boundary Value Problems and Elasticity and Their Connection with Failure Initiation: *Springer*; 2012. 185-220 p.

11. Gómez F, Elices M, Planas J. The cohesive crack concept: application to PMMA at −60°C. *Engineering fracture mechanics*. 2005;72(8):1268-85. doi: 10.1016/j.engfracmech.2004.09.005.

12. Elices M, Guinea G, Gomez J, Planas J. The cohesive zone model: advantages, limitations and challenges. *Engineering Fracture Mechanics*. 2002;69(2):137-63. doi: 10.1016/S0013-7944(01)00083-2.

13. Cendón D, Jin N, Liu Y, Berto F, Elices M. Numerical assessment of gray cast iron notched specimens by using a triaxiality-dependent cohesive zone model. *Theoretical and Applied Fracture Mechanics*. 2017;90:259-67. doi: 10.1016/J.TAFMEC.2017.06.001.

14. Gómez F, Elices M, Berto F, Lazzarin P. Fracture of U-notched specimens under mixed mode: experimental results and numerical predictions. *Engineering Fracture Mechanics*. 2009;76(2):236-49. doi: 10.1016/j.engfracmech.2008.10.001.

15. Cendon D, Torabi A, Elices M. Fracture assessment of graphite V‐notched and U‐notched specimens by using the cohesive crack model. *Fatigue & Fracture of Engineering Materials & Structures*. 2015;38(5):563-73. doi: 10.1111/FFE.12264.

16. Blal N, Daridon L, Monerie Y, Pagano S. Artificial compliance inherent to the intrinsic cohesive zone models: criteria and application to planar meshes. *International Journal of Fracture*. 2012;178:71-83. doi: 10.1007/s10704-012-9734-y.

17. Espinosa HD, Zavattieri PD. A grain level model for the study of failure initiation and evolution in polycrystalline brittle materials. Part I: Theory and numerical implementation. *Mechanics of Materials*. 2003;35(3-6):333-64. doi: 10.1016/S0167-6636(02)00285-5.

18. Foulk III J. An examination of stability in cohesive zone modeling. *Computer Methods in Applied Mechanics and Engineering*. 2010;199(9-12):465-70. doi: 10.1016/j.cma.2009.08.025.

19. Akbardoost J, Rastin A. Scaling effect on the mixed-mode fracture path of rock materials. *Physical Mesomechanics*. 2016;19(4):441-51. doi: 10.1134/S102995991604010X.

20. Aliha MRM, Ayatollahi M, Smith D, Pavier M. Geometry and size effects on fracture trajectory in a limestone rock under mixed mode loading. *Engineering Fracture Mechanics*. 2010;77(11):2200-12. doi: 10.1016/J.ENGFRACMECH.2010.03.009.

21. Akbardoost J, Torabi AR, Amani H. Predicting the fracture trajectory in U, VO, and key-hole notched specimens using an incremental approach. *Engineering Fracture Mechanics*. 2018;200:189-207. doi: 10.1016/j.engfracmech.2018.07.012.

22. Hillerborg A, Modéer M, Petersson P-E. Analysis of crack formation and crack growth in concrete

by means of fracture mechanics and finite elements. *Cement and Concrete Research*. 1976;6(6):773-81. doi: 10.1016/0008-8846(76)90007-7.

23. Bažant Z, Pfeiffer P. Shear fracture tests of concrete. *Materials and Structures*. 1986;19(2):111-21. doi: 10.1007/BF02481755.

24. Carpinteri A. Decrease of apparent tensile and bending strength with specimen size: two different explanations based on fracture mechanics. *International Journal of Solids and Structures*. 1989;25(4):407-29. doi: 10.1016/0020-7683(89)90056-5.

25. Kullmer G, Richard H. Influence of the root radius of crack-like notches on the fracture load of brittle components. *Archive of Applied Mechanics*. 2006;76(11):711-23. doi: 10.1007/s00419-006-0089-6.

26. Moës N, Dolbow J, Belytschko T. A finite element method for crack growth without remeshing. *International Journal for Numerical Methods in Engineering*. 1999;46(1):131-50. doi: 10.1002/(SICI)1097-0207(19990910)46:1<131::AID-NME726>3.0.

27. Dolbow J, Moës N, Belytschko T. An extended finite element method for modeling crack growth with frictional contact. *Computer Methods in Applied Mechanics and Engineering*. 2001;190(51-52):6825-46. doi: 10.1016/S0045-7825(01)00260-2.

28. Areias PMA, Belytschko T. Analysis of three-dimensional crack initiation and propagation using the extended finite element method. *International Journal for Numerical Methods in Engineering*. 2005;63(5):760-88. doi: 10.1002/nme.1305.

29. Belytschko T, Black T. Elastic crack growth in finite elements with minimal remeshing. *International Journal for Numerical Methods in Engineering*. 1999;45(5):601-20. doi: 10.1002/(SICI)1097-0207(19990620)45:5<601::AID-NME598>3.0.CO;2-S.

30. Bahrami B, Ayatollahi MR, Torabi A. Predictions of fracture load, crack initiation angle, and trajectory for V-notched Brazilian disk specimens under mixed mode I/II loading with negative mode I contributions. *International Journal of Damage Mechanics*. 2018;27(8):1173-91. doi: 10.1177/1056789517726360.

31. Dong L, Yang T, Wang K, Atluri SN. A new Fragile Points Method (FPM) in computational mechanics, based on the concepts of Point Stiffnesses and Numerical Flux Corrections. *Engineering Analysis with Boundary Elements*. 2019;107:124-33.

32. Yang T, Dong L, Atluri SN. A simple Galerkin meshless method, the Fragile Points method using point stiffness matrices, for 2D linear elastic problems in complex domains with crack and rupture propagation. *International Journal for Numerical Methods in Engineering*. 2021;122(2):348-85. doi: 10.1002/nme.6540.

33. Atluri SN, Zhu T. A new meshless local Petrov-Galerkin (MLPG) approach in computational mechanics. *Computational Mechanics*. 1998;22(2):117-27. doi: 10.1007/S004660050346.

34. Belytschko T, Lu YY, Gu L. Element-free Galerkin methods. *International Journal for Numerical Methods in Engineering*. 1994;37(2):229-56. doi: 10.1002/nme.1620370205.

35. Guan Y, Grujicic R, Wang X, Dong L, Atluri SN. A new meshless "fragile points method" and a local variational iteration method for general transient heat conduction in anisotropic nonhomogeneous media. Part II: Validation and discussion. *Numerical Heat Transfer, Part B: Fundamentals*. 2020;78(2):86-109. doi: 10.1080/10407790.2020.1747283.

36. Guan Y, Grujicic R, Wang X, Dong L, Atluri SN. A new meshless "fragile points method" and a local variational iteration method for general transient heat conduction in anisotropic nonhomogeneous media. Part I: Theory and implementation. *Numerical Heat Transfer, Part B: Fundamentals*.

2020;78(2):71-85. doi: 10.1080/10407790.2020.1747278.

37. Guan Y, Dong L, Atluri SN. A new meshless Fragile Points Method (FPM) with minimum unknowns at each point, for flexoelectric analysis under two theories with crack propagation, I: Theory and implementation. *Journal of Mechanics of Materials and Structures*. 2021;16(2):159-95. doi: 10.2140/jomms.2021.16.197.

38. Rivière B. Discontinuous Galerkin methods for solving elliptic and parabolic equations: theory and implementation: *SIAM*; 2008.

39. Abedi R, Hawker MA, Haber RB, Matouš K. An adaptive spacetime discontinuous Galerkin method for cohesive models of elastodynamic fracture. *International Journal for Numerical Methods in Engineering*. 2010;81(10):1207-41. doi: 10.1002/NME.2723.

40. Nguyen VP. Discontinuous Galerkin/extrinsic cohesive zone modeling: Implementation caveats and applications in computational fracture mechanics. *Engineering Fracture Mechanics*. 2014;128:37-68. doi: 10.1016/j.engfracmech.2014.07.003.

41. Radovitzky R, Seagraves A, Tupek M, Noels L. A scalable 3D fracture and fragmentation algorithm based on a hybrid, discontinuous Galerkin, cohesive element method. *Computer Methods in Applied Mechanics and Engineering*. 2011;200(1-4):326-44. doi: 10.1016/j.cma.2010.08.014.

42. Simo JC, Ju J. Strain-and stress-based continuum damage models—I. Formulation. *International Journal of Solids and Structures*. 1987;23(7):821-40. doi: 10.1016/0895-7177(89)90117-9.

43. Krajcinovic D, Fonseka G. The continuous damage theory of brittle materials, part 1: general theory. 1981. doi: 10.1115/1.3157739.

44. Marsden JE, Hughes TJ. Mathematical foundations of elasticity: *Courier Corporation*; 1994.

45. Versino D, Mourad HM, Dávila CG, Addessio FL. A thermodynamically consistent discontinuous Galerkin formulation for interface separation. *Composite Structures*. 2015;133:595-606. doi: 10.1016/j.compstruct.2015.07.080.

46. ABAQUS 6.12. Abaqus Analysis User's Manual, Dassault Systèmes Simulia Corp., Providence, RI, USA; 2012.

47. Tsuji, Iwase, Ando. An investigation into the location of crack initiation sites in alumina, polycarbonate and mild steel. *Fatigue & Fracture of Engineering Materials & Structures*. 1999;22(6):509-17. doi: 10.1046/j.1460-2695.1999.00181.x.

48. Gogotsi GA. Fracture toughness of ceramics and ceramic composites. *Ceramics International*. 2003;29(7):777-84. doi: 10.1016/S0272-8842(02)00230-4.

49. Torabi AR, Fakoor M, Pirhadi E. Fracture analysis of U-notched disc-type graphite specimens under mixed mode loading. *International Journal of Solids and Structures*. 2014;51(6):1287-98. doi: 10.1016/j.ijsolstr.2013.12.024.

50. Alfano G, Crisfield MA. Finite element interface models for the delamination analysis of laminated composites: mechanical and computational issues. *International Journal for Numerical Methods in Engineering*. 2001;50(7):1701-36. doi: 10.1002/nme.93.